# Third-order scale-independent WENO-Z scheme achieving optimal order at critical points

Qin Li, Xiao Huang, Pan Yan, Yi Duan

School of Aerospace Engineering, Xiamen University, Xiamen, Fujian, 361102, China

**Abstract**: As we found previously, when critical points occur within grid intervals, the accuracy relations of smoothness indicators of WENO-JS would differ from that assuming critical points occurring on grid nodes, and accordingly the global smoothness indicator in WENO-Z scheme will differs from the original one. Based on above understandings, we first discuss several issues regarding current third-order WENO-Z improvements (e.g. WENO-NP3, -F3, -NN3, -PZ3 and –P+3), i.e. the numerical results with scale dependency, the validity of analysis assuming critical points occurring on nodes, and the sensitivity regarding computational time step and initial condition in order convergence studies. By analyses and numerical validations, the defections of present improvements are demonstrated, either scale-dependency of results or failure to recover optimal order when critical points occurring at half nodes, and then a generic analysis is provided which considers the first-order critical point occurring within grid intervals. Based on achieved theoretical outcomes, two scale-independent, third-order WENO-Z schemes are proposed which could truly recover the optimal order at critical points: the first one is acquired by limited expansion of grid stencil, deriving new global smoothness indicator and incorporating with the mapping function; the second one is achieved by further expanding grid stencil and employing a different global smoothness indicator. For validating, typical 1-D scalar advection problems, 1-D and 2-D problems by Euler equations are chosen and tested. The consequences verify the optimal order recovery at critical points by proposed schemes and show that: the first scheme outperforms aforementioned third-order WENO-Z improvements in terms of numerical resolution, while the second scheme indicates weak numerical robustness in spite of improved resolution and is mainly of theoretical significance.

## 1 Introduction

It is well known that the weighted essentially non-oscillatory (WENO) scheme [1, 2] is one of the most popular methods to improve the order of total variation diminishing (TVD) scheme. Especially, the implementation of WENO-JS [2] has been widely practiced by applications and underwent constant development. Among the developments, the mapping method [3] (with the corresponding WENO-M) and WENO-Z approach [4-7] become the two representatives.

In spite of their distinct realizations, WENO-M and WENO-Z originate from the analysis on the order of WENO schemes [3, 4]. Moreover, the accuracy relations of schemes near critical points were particularly concerned, and the inadequate preservation of optimal order there by WENO-JS was manifested. As the remedies, the mapping method employs a function which does not rely on the explicit variable distribution, while WENO-Z scheme devise a global smoothness-indicator $\tau$ which consists of the local smooth indicator $\beta_k^{(r)}$ ($r$ denotes the grid number of candidate scheme) and relates high-order derivatives of variables and use $\tau/\beta_k^{(r)}$ as the component to build the new nonlinear weight. In Refs. [3-4], both methods provide concrete measures for WENO5-JS to recover optimal order in the case of first-order critical points, and notable improvements are achieved on order convergence, numerical resolution, and etc.

Although both methods above deserve great attentions, attentions in this study mainly concern WENO-Z approach. In addition, considering that the numerical efficiency and robustness are particularly concerned by engineering applications, we further focus on the specific third-order development comprising of two candidate schemes on two-point stencil as that of WENO3-JS. The practice to improve aforementioned third-order WENO-Z, especially that by Don et al. [6] which is referred as WENO3-Z subsequently, at least started from Wu et al. [8-9]. They were motivated by the observation that if the derivation of $\tau$ followed the way of WENO5-Z [6], namely $\tau_3 = |\beta_0^{(2)} - \beta_1^{(2)}|$, $\tau_3/\beta_k^{(2)}$ would take the accurate order of $O(\Delta x)$ in smooth region other than $O(\Delta x^2)$ that was required by sufficient condition [4, 8]. Consequently, a new global smoothness-indicator $\tau_N$ [9] was proposed and new $\tau_N^p/\beta_k^{(2)}$ was used as the component to construct the nonlinear weight with $p = \frac{3}{2}$. Next, Xu & Wu [10-11] proposed an analogous component $\tau^*/\beta_k^{(2)p}$ but with $p$ on the denominator, where $\tau^*$ takes $\tau_3$ or $\tau_N$. They further took the idea of Acker et al. [7] to improve practical performance of WENO-Z and proposed a third-order WENO-P+3 based on a slightly modified $\tau_P$ from $\tau_N$. It was indicated [12] that WENO-P+3 cannot recover the third-order in the case of critical points. It is also observed that Gande et al. [13] proposed a $\tau_{F3}$ which is actually very similar to $\tau_N$ and a same form of component $\tau_{F3}^p/\beta_k^{(2)}$ as that in Ref. [9]. There should be other similar works in this regard, however, due to their quite resemblance to the above works, they will not be repeated again.

Although the above improvements which are based on candidate schemes as that of WENO3-JS, the following doubts still exist: (1) All methods reviewed above either employ a dimensional component ($\tau^p/\beta_k^{(2)}$ or $\tau/\beta_k^{(2)p}$) or directly use $\Delta x$. Will corresponding numerical results be scale-dependent when different scales of variable or length are applied? (2) Although above improvements except that of Ref. [12] dedicated to recover the optimal order in the case of first-order critical points, the convergence order based on $L_\infty$ error regarding canonical scalar advection problem was either absent or did not achieve 3 in Refs. [9-11]; although Ref. [13] indicated achieving the order with specific CFL number, it is unknown if the achievement would still succeed if other CFL number is applied. Thorough investigations on optimal order achievement in the occurrence of critical points are inadequate. (3) The current analyses and consequent scheme constructions are carried out by assuming the occurrence of critical points at the reference grid point. Does this assumption strictly hold and would different consequences be possible? Considering the above doubts, further analyses and investigations on WENO3-Z improvement are still worthwhile.

In this paper, we dedicate to develop a scale-independent WENO3-Z improvement with optimal order recovery on first-order critical points. The analysis and development is based on the assumption that the critical points would occur among grid intervals; furthermore, in order to fulfil the target, a recently devised piecewise rational mapping with fine regulation capability [14] is incorporated. The paper is arranged as follows: the consequences of accuracy analysis and requirement of WENO schemes are reviewed in Section 2; several issues regarding improvement of WENO3-Z are analyzed in Section 3; based on the newly obtained understandings, an improved third-order WENO-Z scheme is proposed in Section 4 which incorporates with the mapping in Ref. [14]; in Section 5, numerical examples are provided for the developed scheme; conclusions are drawn in Section 6 at last.

**2 Reviews on accuracy analyses and improvements of WENO3-Z**

In order to facilitate discussion, the formula of WENO-JS [2] is described first. By convention, consider the following scalar, one-dimensional hyperbolic conservation law,

$$u_t + f(u)_x = 0, \tag{1}$$

where $f'(u) > 0$ is assumed. As the grids are discretized as $x_j = j\Delta x$ where $\Delta x$ denotes the interval and $j$ is the grid index, Eq. (1) can be equivalently written as: $(f(u)_x)_j = -\left(h_{j+1/2} - h_{j-1/2}\right)/\Delta x$, where $f(x) = \frac{1}{\Delta x}\int_{x-\Delta x/2}^{x+\Delta x/2} h(x')\,dx'$. Take the conservative scheme $\hat{f}(x)$ as the approximation of $h(x)$, then

$$(f(u)_x)_j \approx -\left(\hat{f}_{j+1/2} - \hat{f}_{j-1/2}\right)/\Delta x \tag{2}$$

Supposing $r$ is the number of grids of candidate scheme, WENO-JS can be formulated as:

$$\hat{f}_{j+1/2} = \sum_{k=0}^{r-1} \omega_k q_k^r \text{ with } q_k^r = \sum_{l=0}^{r-1} a_{kl}^r f(u_{j-r+k+l+1}) \tag{3}$$

where $q_k^r$ is the candidate scheme with $a_{kl}^r$ and $\omega_k$ is the nonlinear weight. $\omega_k$ is derived from the corresponding linear weight $d_k$ as [2]

$$\omega_k = \alpha_k / \sum_{l=0}^{r-1} \alpha_l \tag{4}$$

where $\alpha_k$ is regarded the non-normalized weight, $\varepsilon = 10^{-6} \sim 10^{-7}$ for WENO-JS and $\beta_k^{(r)}$ is the smoothness indicator. $\beta_k^{(r)}$ can be formulated in positive semi-definite quadratic form as [2]:

$$\beta_k^{(r)} = \sum_{m=0}^{r-2} c_m^r \left(\sum_{l=0}^{r-1} b_{kml}^r f(u_{j-r+k+l+1})\right)^2 \tag{5}$$

where $c_m^r$ and $b_{kml}^r$ are coefficients. For WENO3-JS and WENO5-JS, $r = 3$ and 5 respectively, and the values of $a_{kl}^r$, $d_k$, $c_m^r$ and $b_{kml}^r$ are tabulated in Table ? and ? in Appendix for reference.

For WENO-JS, $\alpha_k$ take the form as $\alpha_k = d_k/\left(\varepsilon + \beta_k^{(r)}\right)^2$ [2]; for WENO-Z, $\alpha_k$ would take the form such as

$$\alpha_k = d_k \left(1 + \left(\tau/\beta_k^{(r)}\right)^p\right) \tag{6}$$

in WENO5-Z [4] where $\tau$ is sometimes referred as the global smoothness-indicator and $p = 1$ or 2. $\tau$ is so designed that it is generally much smaller than $\beta_k^{(r)}$ by several orders, and therefore $\alpha_k$ would be more close to $\alpha_k$ to yield the more accuracy of WENO-Z. Besides, the implementation of WENO-Z provides possibility to preserve optimal order at critical points.

(1) Accuracy analysis and requirements of WENO

Although accuracy relations had been provided in Ref. [2], they were re-visited in Refs. [3-4], and the following conditions were acquired to make the scheme achieve the design order (2$r$-1).

Necessary and sufficient conditions:

$$\begin{cases} \sum_{k=0}^{r-1} A_k(\omega_k^+ - \omega_k^-) = O(\Delta x^r) \\ \omega_k^\pm - d_k = O(\Delta x^{r-1}) \end{cases} \tag{7}$$

where the superscript "$\pm$" denotes the location $x_{j\pm\frac{1}{2}}$.

Sufficient condition:
$$\omega_k^\pm - d_k = O(\Delta x^r) \tag{8}$$

It is worth mentioning the concrete form of $\omega_k$ is not required at this moment. As indicated in Ref. [4], when derivatives of $f$ exist and are continuous, the validity of Eq. (7) can be testified for WENO-JS through symbolic derivation. Regarding the third-order WENO3-Z [6] where Eq. (6) is used with $p=1$ and $\tau = \tau_3$ with $\tau_3$ shown in Eq. (9), Wu & Zhao [8] showed that without the presence of critical points, Eq. (8) was not satisfied, whereas Xu & Wu [9] pointed out Eq. (7) still held or the third-order was established yet. However, all references [9-13] agreed that the achievement of optimal order would be violated when critical points occur with certain order $n$. From the references up-to-date, the analyses were carried out by assuming $f^{(i)}(x_j) = 0$ for $i \le n$ under the scenario of Eq. (2). It can be seen that corresponding improvements to recover the optimal order usually fall in two aspects: the first one aims to satisfy Eq. (7) [10-11] and the second one aims to satisfy Eq. (8) [9,13]; and the concrete realizations include: the new formulation of $\tau$ and integration of $\tau$ and $\beta_k^{(r)}$ in $\alpha_k$.

(2) Typical WENO3-Z improvements

As mentioned in the introduction, we only concern the WENO3-Z improvements comprising of candidate schemes with two grids, and reviews are merely made on those referred previously. First, the formulation of global smoothness-indicator $\tau$ referred can be summarized as:

$$\begin{cases} \tau_3 = \left|\beta_0^{(2)} - \beta_1^{(2)}\right|^{[5]}, & \tau_N = \left|\frac{1}{2}\left(\beta_0^{(2)} + \beta_1^{(2)}\right) - \beta_1^{(3)}\right|^{[8-9]} \\ \tau_{F3} = \left|\frac{1}{2}\left(\beta_0^{(2)} + \beta_1^{(2)}\right) - \dot\beta_1^{(3)}\right|^{[13]}, & \tau_p = \left|\frac{1}{2}\left(\beta_0^{(2)} + \beta_1^{(2)}\right) - \frac{1}{4}(f_{i-1} - f_{i+1})^2\right|^{[12]} \end{cases} \tag{9}$$

where $\dot\beta_1^{(3)}$ is one similar to $\beta_1^{(3)}$ as $\dot\beta_3 = \frac{1}{4}(f_{i-1} - f_{i+1})^2 + \frac{1}{12}(f_{i-1} - 2f_i + f_{i+1})^2$. The motivation of various $\tau$ is to make it have higher orders of $\Delta x$ as possible than that of $\beta_k^{(2)}$, and therefore the satisfaction of Eq. (7) or (8) would possible. In Ref. [12], Xu & Wu pointed out that $\tau_N$ and $\tau_p$ are actually $\frac{10}{12}(f_{i-1} - 2f_i + f_{i+1})^2$ and $\frac{3}{12}(f_{i-1} - 2f_i + f_{i+1})^2$ respectively. Observing this indication, we further investigate $\tau$ starting from the following general formulation:

$$\tau^* = \left|a_1(f_j - f_{j-1})^2 + a_2(f_{j+1} - f_j)^2 + a_3(f_{j+1} - f_{j-1})^2 + b(f_{j+1} - 2f_j + f_{j-1})^2\right| \tag{10}$$

which includes all discretizations of possible derivatives at the stencil $(x_{j-1}, x_j, x_{j+1})$ in quadratic form. Then, the coefficients $a_k$ and $b$ can be defined to make $\tau^*$ achieve desired order.

(a) $\tau^*$ with the error as $O(\Delta x^3)$

It is easy to find that $a_3$ should satisfy: $a_3 = -\frac{1}{4}(a_1 + a_2)$. If $-a_1 = a_2 = a$, then $a_3 = 0$ and $\tau^*$ can be derived as $\tau^* = \left|a(f_{j+1} - f_{j-1})(f_{j+1} - 2f_j + f_{j-1}) + b(f_{j+1} - 2f_j + f_{j-1})^2\right|$; furthermore when $b = 0$, $\tau^* = \left|a(f_{j+1} - f_{j-1})(f_{j+1} - 2f_j + f_{j-1})\right|$ which is similar to $\tau_3$. Besides, one can check that when $b = -a$, $\tau^* = \left|a(f_j - f_{j-1})(f_{j+1} - 2f_j + f_{j-1})\right|$, and when $b = a$, $\tau^* = \left|a(f_{j+1} - f_j)(f_{j+1} - 2f_j + f_{j-1})\right|$.

(b) $\tau^*$ with the error as $O(\Delta x^4)$

It is easy to find that $a_k$ should satisfy: $a_1 = a_2 = a$ and $a_3 = -\frac{1}{2}a$, then $\tau^*$ becomes:

$$\tau^* = \left|c(f_{j+1} - 2f_j + f_{j-1})^2\right| \text{ with } c = \frac{1}{2}(a + 2b). \text{ One can check } \{\tau_N, \tau_p, \tau_{F3}\} \text{ in Eq. (8) can be}$$

re-produced in this form under $c = \{\frac{10}{12}, \frac{3}{12}, \frac{2}{12}\}$, and the first two of which agrees with that reported in Ref. [12]. It is worth mentioning that when $a = -2b$, a trivial solution $\tau^* = 0$ will be obtained.

Hence, one can observe $\tau^*$ and that in Eq. (8) consists of the multiplication of discretizations of derivatives, either that of first and second derivatives or the second derivative square.

Next, the integrations of $\tau$ and $\beta_k^{(r)}$ will be reviewed, through which the non-normalized weight $\alpha$ is defined. Besides that of WENO3-Z, i.e. Eq. (6) at $p = 1$, two other types of implementations were observed as

$$\begin{cases} \alpha_k = d_k \left(1 + \tau^p / \beta_k^{(2)}\right) & \text{[8-9,13]} \quad (11.a) \\ \alpha_k = d_k \left(1 + \tau / \beta_k^{(2)p}\right) & \text{[10-11]} \quad (11.b) \end{cases}$$

through which aforementioned Eq. (7) or (8) was said to be satisfied at first-order critical points. Moreover, Ref. [13] even referred that Eq. (10.a) could achieve the third-order accuracy at second-order critical points providing $\tau = \tau_{F3}$ and $p = 3/2$. For reference, a short summary of characteristics of WENO3-Z and its improvements with progress on critical points are given in Table 1, where $r_{CP,MAX}$ denotes the maximum order of critical points for the scheme to recover the optimal order 3 reported in the references, "$p$ by Ref." means $p$ is derived by corresponding references, and "$p$ by ACP$_\lambda$" means $p$ is derived by new analysis based on critical points occurring within intervals discussed later.

Table 1 Short summary of WENO3-Z and improvements with achievement on critical points

| Scheme | $\beta_k^{(r)}$ | $\alpha_k$ | $p$ by Ref. | $p$ by ACP$_\lambda$ | $r_{CP,MAX}$ |
|---|---|---|---|---|---|
| WENO3-Z [6] | $\tau_3$ | Eq. (6) | $p = 1$ | N/A | 0 |
| WENO-NP3[9] | $\tau_N$ | Eq. (11.a) | $p = 3/2$ | $p \geq 3/2$ | 2 |
| WENO-F3[13] | $\tau_{F3}$ | Eq. (11.a) | $p = 3/2$ | $p \geq 3/2$ | 2 |
| WENO-NN3[10] | $\tau_N$ | Eq. (11.b) | $p \leq 3/4$ | $p \leq 1/2$ | 1 |
| WENO-PZ3[11] | $\tau_3$ | Eq. (11.b) | $p \leq 1/2$ | $p \leq 1/2$ | 1 |

In order to facilitate subsequent discussions, $\alpha_k$ of WENO-P+3 [12] is also shown here as: $\alpha_k = d_k \left(1 + \frac{\tau_p}{\beta_k + \varepsilon} + \lambda_{p+3} \frac{\beta_k + \varepsilon}{\tau_p}\right)$ where $\lambda_{p+3} = \Delta x^{1/6}$ and $\varepsilon = 10^{-40}$. As indicated in the reference, WENO-P+3 dedicates to the improvement of numerical resolution, not the recovery of optimal order at critical points.

(3) Propositions regarding ratio of nonlinear weights on smooth and less smooth stencils

As indicated in Refs. [4, 8-13], the improvement on structural resolution is quite concerned. In order to shed light in this regard, the parametric study of the ratio of nonlinear weights on smooth and less smooth stencils is made in Refs. [10-12] on ones such as $p$ in Eq. (11), through which the influence of parameter on structural resolution is acquired. The following propositions have been proposed or indicated.

Suppose $S_C$ and $S_D$ are two substencils of the same pattern, where the subscripts "C, D" indicates the variable distribution at $S_C$ is smoother than at $S_D$, or $\beta_{k,C}^{(r)} < \beta_{k,D}^{(r)}$, and suppose $\alpha_k$ is the non-normalized weight while $\omega_k$ is the normalized weight.

**Proposition 1** [11]. Consider $\alpha_k$ with the form $\alpha_k = d_k\left(1 + \tau/\beta_k^{(r)p_i}\right)$, then $\left[\frac{\omega_{k,D}}{\omega_{k,C}}\right]_{p_1} > \left[\frac{\omega_{k,D}}{\omega_{k,C}}\right]_{p_2}$ for $\tau > 0$, $0 < p_1 < p_2 \leq 1$ and $0 < \beta_{k,C}^{(r)} < \beta_{k,D}^{(r)} \leq e$.

In Ref. [11], the proposition was given in the case of $p_2=1$ but in the absence of condition "$0 < \beta_{k,C}^{(r)} < \beta_{k,D}^{(r)} \leq e$". The proof of the proposition is given in Appendix II.

**Proposition 2.** Consider $\alpha_k$ with the form $\alpha_k = d_k\left(1 + \tau^{p_i}/\beta_k^{(r)}\right)$, then $\left[\frac{\omega_{k,D}}{\omega_{k,C}}\right]_{p_1} < \left[\frac{\omega_{k,D}}{\omega_{k,C}}\right]_{p_2}$ for $1 \leq p_1 < p_2$ and $0 < \tau < 1$; $\left[\frac{\omega_{k,D}}{\omega_{k,C}}\right]_{p_1} > \left[\frac{\omega_{k,D}}{\omega_{k,C}}\right]_{p_2}$ for $1 \leq p_1 < p_2$ and $1 < \tau$; $\left[\frac{\omega_{k,D}}{\omega_{k,C}}\right]_{p_1} = \left[\frac{\omega_{k,D}}{\omega_{k,C}}\right]_{p_2}$ for $1 \leq p_1 < p_2$ and $\tau = 1$.

This proposition was not explicitly stated in Refs. [9, 13], however, it would be regarded as the counterpart of Proposition 1. Its proof is provided in Appendix II for completeness.

**Proposition 3** [12]. Consider $\alpha_k$ with the form $\alpha_k = d_k\left(1 + c_i \times \tau/\beta_k^{(r)}\right)$, then $\left[\frac{\omega_{k,D}}{\omega_{k,C}}\right]_{c_1} < \left[\frac{\omega_{k,D}}{\omega_{k,C}}\right]_{c_2}$ for $c_1 > c_2$.

This proposition was mentioned in Ref. [12] in similar manner, whose proof is shown in Appendix II again for completeness.

Borges et al. [4] once mentioned that when $p$ in Eq. (6) takes the larger 2 than 1, the relative importance of discontinuous stencil is decreased and corresponding scheme appears more dissipative. As a correspondence, the following proposition is presented with the proof shown in Appendix II as well.

**Proposition 4.** Consider $\alpha_k$ with the form $\alpha_k = d_k\left(1 + \left(\tau/\beta_k^{(r)}\right)^{p_i}\right)$, then $\left[\frac{\omega_{k,D}}{\omega_{k,C}}\right]_{p_1} < \left[\frac{\omega_{k,D}}{\omega_{k,C}}\right]_{p_2}$ for $p_1 > p_2$ providing $0.278 < \frac{\tau}{\beta_{k,D}^{(r)}} < \frac{\tau}{\beta_{k,C}^{(r)}}$.

**3 Issues of WENO3-Z improvements and new analysis based on critical points within intervals**

In spite of progress regarding WENO3-Z improvements, we find that some doubts and unexpected results different from that in the original literatures exist after careful investigations. To this end, discussions are made on these issues at first.

3.1 Issues to be clarified regarding current WENO3-Z improvements

(1) Scale-independency of scheme

Considering Eq. (11), one will find $\tau^p/\beta_k^{(r)}$ is dimensional providing $p \neq 1$, which means $\alpha_k$ would have inconsistent values when different variable scales of $f$ are used to solve a same problem, and essentially different numerical results will be engendered. To clarify this issue, we take Shu-Osher problem in Section 5 as illustration. Suppose a new variable scale is employed such that a ratio of new and original initial state variables holds as: $r_{var} = \left(\frac{\rho_{new}}{\rho}\right)_0 = \left(\frac{p_{new}}{p}\right)_0$ where the initial state is indicated by the subscript "0" and the original state means that described in Section 5; meanwhile, the initial velocity remains invariable. For brevity, the computation and corresponding result of the original state are generally referred as the original one. We refer a scheme as variable scale-independent in view of that, when the new computation is accomplished by using the same time step $\Delta t$ as the original one, the result of a scale-independent scheme should have the same density and pressure distributions as that of the original initial conditions as long as they are scaled again by multiplying $1/r_{var}$. In purpose of demonstration, three schemes are tested: WENO3-JS, WENO-P+3 and WENO-F3. As just discussed, theoretically, the two formers are variable scale-independent and the last one is variable scale-dependent. The computations are carried out under $\Delta t = 0.0015$ and $r_{var} = 1/10$ on 400 grids. The rescaled new result compared with the original one of WENO-P+3 is shown in Fig. 1(a), while those regarding WENO-F3 is shown in Fig. 1(b). One can see that the new result of WENO-P+3 coincides with the original one whereas those of WENO-F3 indicate obvious difference. The performance of WENO3-JS is the same as that of WENO-P+3 and is omitted thereby. Hence, the variable scale-independence of WENO3-JS and WENO-P+3 and scale-dependence of WENO-F3 are demonstrated.

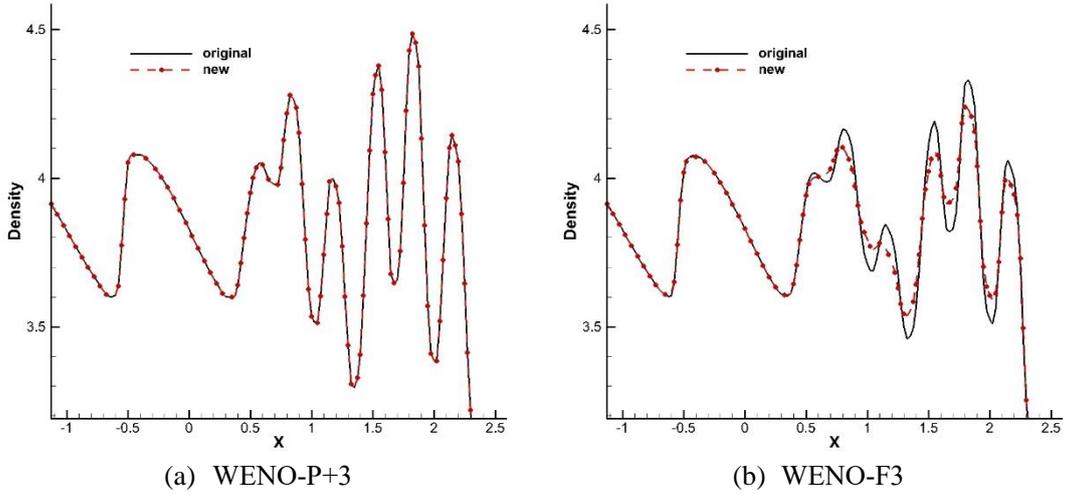

(a) WENO-P+3  (b) WENO-F3

Fig. 1 Local view of density distributions of Shu-Osher problem of WENO3-Z improvements by respectively employing the new and original initial states (Grids: 400; $\Delta t = 0.0015$)

Next, we consider another kind of scale-independence, namely length scale-independence. Take aforementioned one-dimensional problem for illustration. If $x$ coordinate is scaled by a factor $r_{\Delta x}$ as $x_{new} = r_{\Delta x} \cdot x$, we refer a scheme as length scale-independent in view of that, when the new computation is accomplished by using the same CFL number for $\Delta t$ as the original one, the result of a length scale-independent scheme should have the same distributions as that of the original one as long as the coordinates are scaled back by multiplying $1/r_{\Delta x}$. The original computation is the same one as above, and the new computation is carried out on coordinates at $r_{\Delta x} = 100$. The new results of WENO-P+3 is shown in Fig. 2(a) with the comparison with the original one, while those of WENO3-JS is shown in Fig. 2(b). The figure tells that the new result of WENO-P+3 disagrees with the original one (due to direct use of $\lambda_{p+3} = \Delta x^{1/6}$ in the scheme), whereas those of WENO3-JS coincide with each other. The performance of WENO-F3 is similar to that of WENO3-JS and is omitted as well. Hence, the length scale-independence of WENO3-JS and WENO-F3 and scale-dependence of WENO-P+3 are manifested.

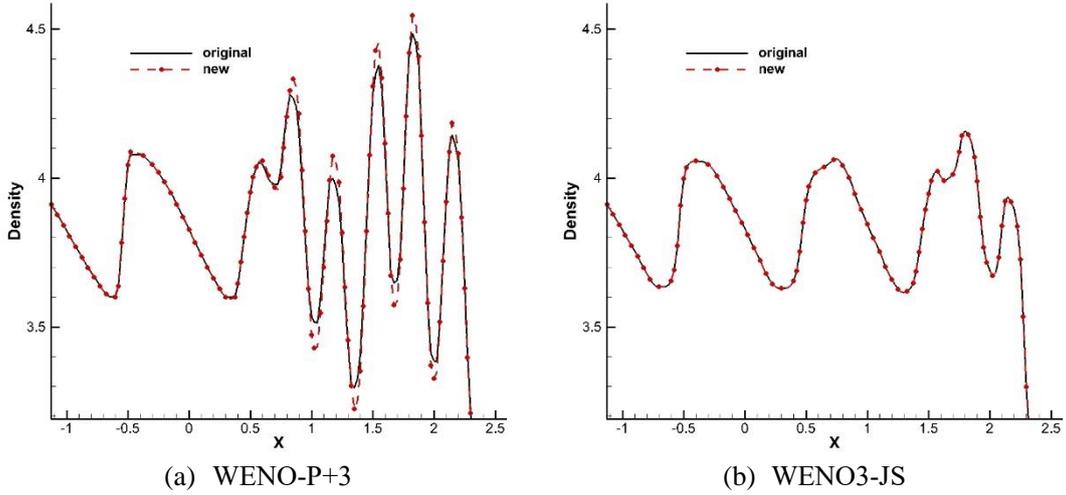

(a) WENO-P+3  (b) WENO3-JS

Fig. 2 The density distributions of Shu-Osher problem of WENO3-Z improvements by using the new and original coordinates respectively (Grids: 400; $CFL = 0.06$)

Similarly, one can check other WENO3-Z improvements in Table 1 except WENO3-Z are all variable scale-dependent. For convenience, we refer the scale-independent property, either variable or length ones, with the abbreviation as SCI.

(2) Analysis and corresponding recipes to recover accuracy at critical points

To facilitate discussion, the scenario on which Eqns. (2) -(5) are introduced is employed. At least in Ref. [9], Wu et al. found that when $f_j' = 0$, the proposed $\tau_N$ would yield: $\tau_N/\beta_k^{(r)} \sim O(\Delta x^4)/O(\Delta x^4)$ which nullified the attempt to satisfy Eq. (8). As the remedy, they employed the recipe as Eq. (11.a) and gave a solution of $p$ as 1/2. One can verify that $\alpha_k = d_k(1 + O(\Delta x^2))$ thereafter and Eq. (8) be acquired. As shown in Ref. [9], the achievement of third-order in the case of certain first-order critical points was indicated. It is worth mentioning that in their analysis, the key point is to assume the critical point occurs at the node, namely $x_j$ in Eq. (2). This assumption is also widely adopted by other investigations [10-13].

If the analysis by Wu et al. [9] is valid, the analogous measure should be applicable for WENO3-Z where $\alpha_k$ is derived by Eq. (6) at $p = 1$ and $\tau = \tau_3$ in Eq. (9). One can easily find that when non-critical points occur, $\tau_3/\beta_k^{(2)} \sim O(\Delta x)$; when the first-order critical point emerges ($f_j' = 0$), $\tau_3/\beta_k^{(2)} \sim O(\Delta x)$. Hence, a straightforward recipe for WENO3-Z is to still employ Eq. (6) but having $p = 2$, through which $\alpha_k = d_k(1 + O(\Delta x^2))$ and Eq. (8) be satisfied, and therefore the third-order is expected to recover. Such practice is easy to numerically test, and canonical 1-D scalar advection problem with the initial condition $u(x,0) = \sin(\pi x - \frac{\sin(\pi x)}{\pi})$ is usually adopted. In the case the first-order critical point exists where its second- and third-order derivatives are non-zero, and the details of which are given in Section 5. The results of errors (especially $L_\infty$ error) and convergence of numerical order are shown in Table 2. It is interesting to note that aforementioned trial of WENO3-Z based on assumption critical point occurring at nodes cannot achieve the optimal order. As will be shown in Section 3.2, another practice to improve WENO3-Z is proposed by employing Eq. (11.a) and $p = 3/2$. Corresponding results are also shown in Table 2, which indicates the achievement of optimal order. Hence aforementioned analysis based on critical point at nodes seems to be questionable.

Table 2 Numerical errors and order convergence of WENO3-Z and Eq. (11.a) at $p = 3/2$ and

$\tau = \tau_3$ by 1-D scalar advection problem with $u(x, 0) = \sin(\pi x - \sin(\pi x)/\pi)$ and CFL=0.4

| N | Δt | WENO3-Z | | Eq. (11.a) at $p = 3/2$ & $\tau = \tau_3$ | |
|---|---|---|---|---|---|
| | | $L_1$-error | $L_1$-order | $L_1$-error | $L_1$-order |
| 10 | 0.08 | 2.1044E-01 | | 1.9317E-01 | |
| 20 | 0.04 | 7.4181E-02 | 1.504 | 4.6739E-02 | 2.047 |
| 40 | 0.02 | 2.3950E-02 | 1.631 | 5.4225E-03 | 3.107 |
| 80 | 0.01 | 5.7003E-03 | 2.070 | 5.3359E-04 | 3.345 |
| 160 | 0.005 | 1.2774E-03 | 2.157 | 6.0655E-05 | 3.137 |
| 320 | 0.0025 | 2.7330E-04 | 2.224 | 8.1658E-06 | 2.892 |
| 640 | 0.00125 | 5.7674E-05 | 2.244 | 8.1066E-07 | 3.332 |
| | | $L_\infty$-error | $L_\infty$-order | $L_\infty$-error | $L_\infty$-order |
| 10 | 0.08 | 4.6613E-01 | | 4.2350E-01 | |
| 20 | 0.04 | 1.8677E-01 | 1.319 | 1.3421E-01 | 1.657 |
| 40 | 0.02 | 7.2992E-02 | 1.355 | 2.7125E-02 | 2.306 |
| 80 | 0.01 | 2.7400E-02 | 1.413 | 3.3319E-03 | 3.025 |
| 160 | 0.005 | 9.9731E-03 | 1.458 | 4.5245E-04 | 2.880 |
| 320 | 0.0025 | 3.5506E-03 | 1.489 | 9.2707E-05 | 2.287 |
| 640 | 0.00125 | 1.2759E-03 | 1.476 | 4.2499E-06 | 4.447 |

Considering that the formulations of Wu et al. [9] and Gande et al. [13] take Eq. (11.a), in order to further explore the validity of analysis based on critical point at nodes, we investigate the practice of WENO-NN3 where Eq. (11.b) is employed. Still assuming the occurrence of critical point at $x_j$, Xu & Wu [10] mentioned WENO-NN3 would achieve the optimal third-order at the first-order critical points with $f_j'' \& f_j''' \neq 0$ providing $p \leq 3/4$. As indicated in the introduction, current WENO3-Z improvements including WENO-NN3 favors order convergence study based on $L_1$ or $L_2$ but other than $L_\infty$, however, results from the latter should be more plausible because the occurrence of critical points in tests is usually as seldom as 1~2 time(s). In above consideration, we test the order convergence of WENO-NN3 at $p = 3/4$ using the same problem as above, and the results are shown in Table 3. It is definite that WENO-NN3 at $p = 3/4$ cannot achieve the third-order. As just mentioned, we propose a new analysis in Section 3.2, and according to consequence of the analysis, the request for WENO-NN3 to recover the optimal order should be $p \leq 1/2$. Accordingly, we make the same tests just now with results shown in Table 3 also. The consequence tells that WENO-NN3 can achieve the third-order this time.

Table 3 Numerical errors and order convergence of WENO-NN3 at $p = 3/4$ and $1/2$ by 1-D scalar advection problem with $u(x, 0) = \sin(\pi x - \sin(\pi x)/\pi)$ and CFL=0.4

| N | Δt | WENO-NN3 at $p = 3/4$ | | WENO-NN3 at $p = 1/2$ | |
|---|---|---|---|---|---|
| | | $L_1$-error | $L_1$-order | $L_1$-error | $L_1$-order |
| 10 | 0.08 | 1.6635E-01 | | 1.2411E-01 | |
| 20 | 0.04 | 3.8984E-02 | 2.093 | 2.3192E-02 | 2.419 |
| 40 | 0.02 | 6.1280E-03 | 2.669 | 3.2826E-03 | 2.820 |
| 80 | 0.01 | 8.4603E-04 | 2.856 | 4.0914E-04 | 3.004 |
| 160 | 0.005 | 9.9269E-05 | 3.091 | 5.0711E-05 | 3.012 |

| N | Δt | | | | |
|---|---|---|---|---|---|
| 320 | 0.0025 | 1.2552E-05 | 2.983 | 6.3527E-06 | 2.996 |
| 640 | 0.00125 | 1.3245E-06 | 3.244 | 7.9536E-07 | 2.997 |
| | | $L_\infty$-error | $L_\infty$-order | $L_\infty$-error | $L_\infty$-order |
| 10 | 0.08 | 3.6009E-01 | | 2.5964E-01 | |
| 20 | 0.04 | 1.1862E-01 | 1.602 | 5.6659E-02 | 2.196 |
| 40 | 0.02 | 3.0459E-02 | 1.961 | 7.7109E-03 | 2.877 |
| 80 | 0.01 | 6.0430E-03 | 2.333 | 1.0191E-03 | 2.919 |
| 160 | 0.005 | 1.1405E-03 | 2.405 | 1.2810E-04 | 2.991 |
| 320 | 0.0025 | 1.9695E-04 | 2.533 | 1.6032E-05 | 2.998 |
| 640 | 0.00125 | 3.3734E-05 | 2.545 | 2.0047E-06 | 2.999 |

Hence, we draw the following conclusion: at least for WENO3-Z improvements, the commonly used assumption, namely the first-order critical point only occurring at grid node, is inappropriate and might yield incorrect consequence.

(3) Sensitivity of time step and initial condition in order convergence

In Ref. [13], Gande at al. claimed WENO-F3 would achieve the optimal order at not only the first-order critical point, but also even the second-order critical point. Because WENO-F3 is of little difference with WENO-NP3 as pointed out previously, we take the former as an example to illustrate the sensitivity of time step and initial distribution in order convergence regarding current WENO3-Z improvements.

Still consider the 1-D scalar advection case of $u(x,0) = \sin(\pi x - \sin(\pi x)/\pi)$. A small modification is imposed such that one of its first-order critical point is shifted to the origin: $u(x,0) = \sin(\pi(x + x_c) - \sin(\pi(x + x_c))/\pi)$, where $x_c = 0.5966831869112089637212$. Regarding the new initial distribution, two CFL numbers (or time steps) are employed, i.e. 0.4 and 0.25, and consequent computations for order convergence are made with results shown in Table 4. The table tells that the case at CFL=0.4 shows WENO-F3 achieves a converged order as 3, whereas at the case at CFL=0.25 the scheme does not indicate the third-order convergence.

Table 4 Numerical errors and order convergence of WENO-F3 by 1-D scalar advection problem with $u(x,0) = \sin(\pi(x - x_c) - \sin(\pi(x - x_c))/\pi)$ at CFL=0.4 and 0.25

| N | Δt | CFL=0.4 | | CFL=0.25 | |
|---|---|---|---|---|---|
| | | $L_1$-error | $L_1$-order | $L_1$-error | $L_1$-order |
| 10 | 0.08/0.05 | 1.1736E-01 | | 1.1988E-01 | |
| 20 | 0.04/0.025 | 2.2491E-02 | 2.383 | 2.3027E-02 | 2.380 |
| 40 | 0.02/0.0125 | 3.2078E-03 | 2.809 | 3.3110E-03 | 2.797 |
| 80 | 0.01/0.00625 | 4.1206E-04 | 2.960 | 4.2197E-04 | 2.972 |
| 160 | 0.005/0.003125 | 5.0853E-05 | 3.018 | 5.1989E-05 | 3.020 |
| 320 | 0.0025/0.0015625 | 6.3639E-06 | 2.998 | 6.5009E-06 | 2.999 |
| 640 | 0.00125/0.00078125 | 7.9609E-07 | 2.998 | 8.4067E-07 | 2.951 |
| | | $L_\infty$-error | $L_\infty$-order | $L_\infty$-error | $L_\infty$-order |
| 10 | 0.08/0.05 | 2.4842E-01 | | 2.5309E-01 | |
| 20 | 0.04/0.025 | 5.1731E-02 | 2.263 | 5.2671E-02 | 2.264 |
| 40 | 0.02/0.0125 | 7.1561E-03 | 2.853 | 7.1306E-03 | 2.884 |

| 80 | 0.01/0.00625 | 1.0309E-03 | 2.795 | 1.0222E-03 | 2.802 |
| 160 | 0.005/0.003125 | 1.2822E-04 | 3.007 | 1.6505E-04 | 2.630 |
| 320 | 0.0025/0.0015625 | 1.6035E-05 | 2.999 | 3.0106E-05 | 2.454 |
| 640 | 0.00125/0.00078125 | 2.0047E-06 | 2.999 | 7.3682E-06 | 2.030 |

For further confirmation, WENO-PZ3 is tested using the above case with $p = 1/2$, and corresponding results are shown Table 5. The table tells again that WENO-PZ3 can achieve the optimal order at CFD=0.4 but fails when CFD=0.25

Table 5 Numerical errors and order convergence of WENO-PZ3 by 1-D scalar advection problem with $u(x, 0) = \sin(\pi(x - x_c) - \sin(\pi(x - x_c))/\pi)$ at $CFL$=0.4 and 0.25

| N | $\Delta t$ | CFL=0.4 | | CFL=0.25 | |
|---|---|---|---|---|---|
| | | $L_1$-error | $L_1$-order | $L_1$-error | $L_1$-order |
| 10 | 0.08/0.05 | 1.3725E-01 | | 1.3929E-01 | |
| 20 | 0.04/0.025 | 2.6028E-02 | 2.398 | 2.6133E-02 | 2.414 |
| 40 | 0.02/0.0125 | 3.3886E-03 | 2.941 | 3.4350E-03 | 2.927 |
| 80 | 0.01/0.00625 | 4.2508E-04 | 2.994 | 4.2616E-04 | 3.010 |
| 160 | 0.005/0.003125 | 5.1243E-05 | 3.052 | 5.1591E-05 | 3.046 |
| 320 | 0.0025/0.0015625 | 6.3686E-06 | 3.008 | 6.3676E-06 | 3.018 |
| 640 | 0.00125/0.00078125 | 7.9626E-07 | 2.999 | 7.9672E-07 | 2.998 |
| | | $L_\infty$-error | $L_\infty$-order | $L_\infty$-error | $L_\infty$-order |
| 10 | 0.08/0.05 | 2.8458E-01 | | 2.8688E-01 | |
| 20 | 0.04/0.025 | 6.5547E-02 | 2.118 | 6.4847E-02 | 2.145 |
| 40 | 0.02/0.0125 | 8.0418E-03 | 3.026 | 8.8404E-03 | 2.874 |
| 80 | 0.01/0.00625 | 1.2196E-03 | 2.721 | 1.0941E-03 | 3.014 |
| 160 | 0.005/0.003125 | 1.2698E-04 | 3.263 | 1.4667E-04 | 2.899 |
| 320 | 0.0025/0.0015625 | 1.5994E-05 | 2.989 | 1.9054E-05 | 2.944 |
| 640 | 0.00125/0.00078125 | 2.0034E-06 | 2.996 | 3.4144E-06 | 2.480 |

In fact, one can check that all WENO3-Z improvements in Table 1 would fail to achieve the optimal order in the above case at $CFL = 2^{-n}$ and $n \geq 2$. The reason of the incident is that current investigations are based on analysis assuming critical points occurring at grid node, and therefore they have not realized and cannot avoid consequent theoretical defects. A generic analysis in this regard will be given subsequently.

3.2 New analysis based on critical points occurring within intervals

In Ref. [14], in order to explain why canonical mapping method could not recover the third-order accuracy for WENO3-JS, we first proposed the analysis by considering the occurrence of third-order critical point within the grid interval. In such cases, although the magnitude order in $\Delta x$ of $\beta_k^{(2)}$ remains constant, the concrete formulations of indicators have changed, and accordingly, $\tau_3 = \left|\beta_0^{(2)} - \beta_1^{(2)}\right|$ will be different from that obtained by assuming the critical point happening at grid node. Similar incident happens to $\beta_k^{(3)}$ of WENO5-JS also, and details were provided in the reference. Next, we will follow the idea in Ref. [14] and elaborate a new analysis based on critical point occurring within grid intervals.

Consider the stencil of WENO3-JS, i.e. $\{x_{j-1}, x_j, x_{j+1}\}$ under the scenario of Eqns. (2) -(5). Suppose the critical point occurs somewhere in the stencil whose coordinate can be expressed as:
$$x_c = x_j + \lambda \cdot \Delta x \text{ where } -1 < \lambda < 1. \qquad (12)$$
It is obvious that the occurrence at $x_j$ is included in above consideration by letting $\lambda = 0$. To facilitate discussion, we introduce the notation of $\delta^{(m)_n}$ which regards the discretization of $\left(\frac{\partial^m f}{\partial x^n}\right)_j$ with the nth-order, and $\delta^{(2)_2}$ is: $\delta^{(2)_2} = (f_{j+1} - 2f_j + f_{j-1})$. Using Taylor expansion, the accuracy properties of $\beta_k^{(2)}$ at $\lambda \neq 0$ [14] and $\lambda = 0$, $\tau_3$ and $\left(\delta^{(2)_2}\right)^2$ in the case of critical point are tabulated in Table 6, where "CP$_n$" denotes the n$^{th}$-order critical point. As referred next, "CP$_0$" jcindicates the situation where no critical points occur.

Table 6 Accuracy relations of $\beta_k^{(2)}$, $\tau_3$ and $\left(\delta^{(2)_2}f\right)^2$ in cases of CP$_n$

| | $\lambda$ | CP$_1$ | CP$_2$ |
|---|---|---|---|
| $\beta_0^{(2)}$ | $\neq 0$ | $\frac{1}{4}(2\lambda+1)^2 f''^2_{x_c} \Delta x^4$ $-\frac{1}{6}(2\lambda+1)(3\lambda^2+3\lambda+1)f''_{x_c}f'''_{x_c}\Delta x^5$ $+\left[\begin{array}{c}(\frac{1}{3}\lambda^4+\frac{2}{3}\lambda^3+\frac{7}{12}\lambda^2+\frac{1}{4}\lambda+\frac{1}{24})f''_{x_c}f^{(4)}_{x_c} + \\ (\frac{1}{4}\lambda^4+\frac{1}{2}\lambda^3+\frac{5}{12}\lambda^2+\frac{1}{6}\lambda+\frac{1}{36})f'''^2_{x_c}\end{array}\right]\Delta x^6$ $+O(\Delta x^7)$ | $\frac{1}{36}(3\lambda^2+3\lambda+1)^2 f'''^2_{x_c}\Delta x^6$ $-\frac{1}{72}(3\lambda^2+3\lambda+1)\times$ $(4\lambda^3+6\lambda^2+4\lambda+1)f'''_{x_c}f^{(4)}_{x_c}\Delta x^7$ $+O(\Delta x^8)$ |
| | $= 0$ | $\frac{1}{4}f''^2_j \Delta x^4 - \frac{1}{6}f''_j f'''_j \Delta x^5$ | $\frac{1}{36}f'''^2_j \Delta x^6 - \frac{1}{72}f'''_j f^{(4)}_j \Delta x^7 + O(\Delta x^8)$ |
| $\beta_1^{(2)}$ | $\neq 0$ | $\frac{1}{4}(2\lambda-1)^2 f''^2_{x_c}\Delta x^4$ $-\frac{1}{6}(2\lambda-1)(3\lambda^2-3\lambda+1)f''_{x_c}f'''_{x_c}\Delta x^5$ $+\left[\begin{array}{c}(\frac{1}{3}\lambda^4-\frac{2}{3}\lambda^3+\frac{7}{12}\lambda^2-\frac{1}{4}\lambda+\frac{1}{24})f''_{x_c}f^{(4)}_{x_c} + \\ (\frac{1}{4}\lambda^4-\frac{1}{2}\lambda^3+\frac{5}{12}\lambda^2-\frac{1}{6}\lambda+\frac{1}{36})f'''^2_{x_c}\end{array}\right]\Delta x^6$ $+O(\Delta x^7)$ | $\frac{1}{36}(3\lambda^2-3\lambda+1)^2 f'''^2_{x_c}\Delta x^6$ $-\frac{1}{72}(3\lambda^2-3\lambda+1)\times$ $(4\lambda^3-6\lambda^2+4\lambda-1)f'''_{x_c}f^{(4)}_{x_c}\Delta x^7$ $+O(\Delta x^8)$ |
| | $= 0$ | $\frac{1}{4}f''^2_j\Delta x^4 + \frac{1}{6}f''_j f'''_j \Delta x^5 + O(\Delta x^6)$ | $\frac{1}{36}f'''^2_j \Delta x^6 + \frac{1}{72}f'''_j f^{(4)}_j \Delta x^7 + O(\Delta x^8)$ |
| $\tau_3$ | $\neq 0$ | $-2\lambda f''^2_{x_c}\Delta x^4 +$ $(3\lambda^2+\frac{1}{3})f''_{x_c}f'''_{x_c}\Delta x^5 + O(\Delta x^6)$ | $-\frac{1}{3}(3\lambda^2+1)\lambda f'''^2_{x_c}\Delta x^6 + (\frac{5}{6}\lambda^4 +$ $\frac{7}{12}\lambda^2+\frac{1}{36})\times f'''_{x_c}f^{(4)}_{x_c}\Delta x^7 + O(\Delta x^8)$ |
| | $= 0$ | $\frac{1}{3}f''_j f'''_j \Delta x^5 +$ $(\frac{1}{36}f'''_j f^{(4)}_j + \frac{1}{60}f''_j f^{(5)}_j)\Delta x^7 + O(\Delta x^8)$ | $\frac{1}{36}f'''_j f^{(4)}_j\Delta x^7 + (\frac{1}{1080}f'''_j f^{(6)}_j$ $+\frac{1}{720}f^{(4)}_j f^{(5)}_j)\Delta x^9 + O(\Delta x^{10})$ |
| $\left(\delta^{(2)_2}\right)^2$ | $\neq 0$ | $f''^2_{x_c}\Delta x^4 - 2\lambda f''_{x_c}f^{(3)}_{x_c}\Delta x^5 + O(\Delta x^6)$ | $\lambda^2 f^{(3)2}_{x_c}\Delta x^6 - \frac{1}{6}(6\lambda^2+1)\lambda f^{(3)}_{x_c}f^{(4)}_{x_c}\Delta x^7 +$ $O(\Delta x^8)$ |
| | $= 0$ | $f''^2_j \Delta x^4 + O(\Delta x^6)$ | $\frac{1}{144}f^{(4)2}_j \Delta x^8 + O(\Delta x^{10})$ |

Table 6 indicates:
(1) The order of leading term of $\beta_k^2$ does not change at $\lambda \neq 0$ and $\lambda = 0$ in the occurrence

of $CP_1$, whereas their coefficients are different in two cases.

(2) The order of leading term of $\tau_3$ at $\lambda \neq 0$ differs from that at $\lambda = 0$ in the case of $CP_1$, while such order of $\left(\delta^{(2)_2}\right)^2$ keeps unchanged. As previously mentioned, $\left(\delta^{(2)_2}\right)^2$ corresponds to $\tau_N$, $\tau_{F_3}$ and $\tau_p$.

(3) In order or devise new improvement in the framework of WENO-Z, we should not only consider the accuracy properties at $\lambda = 0$ but also that at $\lambda \neq 0$. For brevity, we abbreviate the analysis considering the occurrence of critical points within intervals as $ACP_\lambda$, while the analysis assuming critical point at node are abbreviated as $ACP_N$ correspondingly.

Based on Table 6 and $ACP_\lambda$, the issues (2) -(3) in Section 3.1 can be clarified as:

(1) According to "(2)" in above, $\alpha_k = d_k(1 + \left(\frac{\tau_3}{\beta_k^{(2)}}\right)^p)$ with $p \geq 2$ would not make Eq. (8) satisfied, and the formulation as $\alpha_k = d_k(1 + \tau_3^p/\beta_k^{(2)})$ would fulfil the job providing $p \geq 3/2$, whose numerical validation has been indicated previously.

(2) Consider the formulation $\alpha_k = d_k(1 + c\left(\delta^{(2)_2}\right)^2/\beta_k^{(2)^p})$ where $c$ is certain coefficient. The solution to fulfill Eq. (8) would be $p \leq 1/2$, which differs from that from the analysis of Ref. [11] based on Eq. (7) and assuming critical points at nodes, namely $p \leq 3/4$, and numerical tests negate the validity of the latter.

(3) Abnormality will arise for current WENO3-Z improvements when $\lambda = \pm 1/2$ with $f''_{x_c}, f'''_{x_c} \neq 0$, in which $\beta_0^{(2)} \sim O(\Delta x^6)$ while $\left(\delta^{(2)_2}\right)^2 \sim O(\Delta x^4)$. Although such occurrence would be seldom and occasional, it is conceivable that the order reduction would occur and the optimal order based on $L_\infty$ error would not reach three.

Base on $ACP_\lambda$ and Table 6, $p$ in Eq. (11) can be derived in terms of satisfying the sufficient condition Eq. (8) and are shown in Table 1. The results and corresponding numerical tests tell that $p$ from the references coincide that from $ACP_\lambda$ except the case of WENO-NN3, and our computation has disapproved the validity of the result by Ref. [11].

## 4 Third-order scale-independent WENO-Z scheme

According to previous discussions, measures in Table 1 have the following inherent deficiencies: scale-dependency and incapability of achieving optimal order when $\lambda = \pm 1/2$ which regards the occurrence of critical points, and the underlying reason lies in the absent understanding of $ACP_\lambda$. In the following, two recipes with SCI achieved are proposed.

4.1 New scale-independent improvement by incorporating mapping

Considering the improved $\tau$ other than $\tau_3$ in Eq. (6), it is natural to think of if there would be other operator such that its dividing with $\beta_k^{(2)}$ would make the sufficient condition (Eq. (8)) hold. If the canonical $\beta_k^{(2)}$ is still employed to measure the smoothness of candidate stencil and the framework of WENO-Z is adopted, the following requirement for $\tau$ should be considered if a SCI scheme is desired:

(1) Observing the use of Eq. (6) with $p = 2$ in the case of WENO5-Z, $\tau/\beta_k^{(2)} \sim O(\Delta x^n)$ would hold under $ACP_\lambda$ and $ACP_N$ with $n \geq 1$.

(2) $\tau$ should have the dimension as $[f]^2$, where "[]" denotes the dimension.

Recalling the discussion of $\tau^*$ in part (2) of Section 2, we first investigate the possible multiplications of two discretizations of derivatives. In this regard, we sort out all derivatives on $x_j$ in stencil $\{x_{j-1}, x_j, x_{j+1}\}$, their discretiztions $\delta^{(m)_n}$ with the highest order and corresponding accuracy relations at CP$_0$ and CP$_1$ (with $\lambda = 0$ and $\neq 0$, see Eq. (12)), where . The derivatives and their discretizations are listed in Eq. (13), and the accuracy relations of multiplications are given in Table 7. The relations of $\delta^{(1)_2}{}^2$ is ignored because of its obvious inadequate accurate order as $\tau$.

$$\begin{cases} \left(\frac{\partial f}{\partial x}\right)_j \approx \frac{1}{2\Delta x}\delta^{(1)_2}: \delta^{(1)_2} = (f_{j+1} - f_{j-1}) \\ \left(\frac{\partial^2 f}{\partial x^2}\right)_j \approx \frac{1}{\Delta x^2}\delta^{(2)_2}: \delta^{(2)_2} = (f_{j+1} - 2f_j + f_{j-1}) \end{cases} \tag{13}$$

Table 7 Accuracy relations of multiplications of two derivative discretizations at CP$_0$ and CP$_1$

| $\tau$ | CP$_0$ | CP$_1$ |
|---|---|---|
| $\delta^{(1)_2} \cdot \delta^{(2)_1}$ | $2f_j'f_j''\Delta x^3 + O(\Delta x^5)$ | $\lambda = 0, \frac{1}{3}f_j''f_j'''\Delta x^5 + O(\Delta x^7);$ |
| | | $\lambda \neq 0, -\frac{1}{2}(2\lambda - 1)f_{j_c}''^2 \Delta x^4 + \frac{3}{2}(\lambda^2 - \lambda + \frac{1}{6})f_{j_c}''f_{j_c}'''\Delta x^5$ |
| $\delta^{(2)_2}{}^2$ | $f_j''^2\Delta x^4 + O(\Delta x^6)$ | $\lambda = 0, f_j''^2\Delta x^4 + O(\Delta x^6);$ |
| | | $\lambda \neq 0, f_{j_c}''^2 \Delta x^4 - 2\lambda f_{j_c}'' f_{j_c}''' \Delta x^5$ |

Considering the accuracy relations of $\beta_k^{(2)}$ in Table 6, one can see that from the above tables there is no possibility to obtain a multiplication which can make at least $\tau/\beta_k^{(2)} \sim O(\Delta x)$ satisfied. Considering the dimension of $\tau$, we further propose the following proposition:

**Proposition 5.** Consider the generic quadratic form of $f$ as $\tau(f) = (f_{j-1}, f_j, f_{j+1})[a_{i_1 i_2}](f_{j-1}, f_j, f_{j+1})^T$ where $[a_{i_1 i_2}]$ is a $3 \times 3$ matrix and $i_1, i_2 = 1\ldots3$ Supposing the first-order critical point would occur at $x_c = x_j + \lambda \cdot \Delta x$ where $-1 < \lambda < 1$ and $\{f_{x_c}' = 0, f_{x_c}'' \& f_{x_c}''' \neq 0\}$, then no solution of $a_{i_1 i_2}$ exists such that Taylor expansion of $\tau(f)$ toward $x_c$ has the leading error as $O(\Delta x^5)$.

The proof of the proposition is shown in Appendix II.

Hence, according to the proposition, there is impossible to derive a SCI, third-order WENO-Z scheme which can recover the optimal order at the first-order critical point. In order to break the predicament, the idea of extending the smoothness indicator is employed, which has once been used in Ref. [14] for mapped WENO3-JS. In detail, only the stencil of original $\beta_1^{(2)}$ is extended to $\{x_{j-1}, x_j, x_{j+1}\}$ and the indicator itself is extended to $\beta_2^{(3)}$ in WENO5-JS, while $\beta_0^{(2)}$ remains unchanged. The benefit of such extension is: (1) The extension of $\beta_1^{(2)}$ indicates the global stencil expands as well, namely $\{x_{j-1}, x_j, x_{j+1}, x_{j+2}\}$, which provides availability to derive a $\tau$ with $\tau/\beta_k^{(2)} \sim O(\Delta x)$ at critical points; (2) The accuracy relation of $\beta_2^{(3)}$ based on ACP$_\lambda$ (see Table 8) indicates the coefficient of the leading term would be non-zero at $\lambda = \pm 1/2$, which mitigates the difficulty to recover optimal order. Although the expansion of stencil indicates more region dependence and computation cost, in practical simulations the fluxes are split and the symmetric scheme is casted for the split one with opposite windward direction, so in this sense the variable on

$y_{j+2}$ still falls in the whole stencil of scheme. Based on above discussion and after careful analysis, we sort out helpful discretizations regarding derivatives at $x_j$ in expended stencil as:

$$\begin{cases} \left(\dfrac{\partial f}{\partial x}\right)_j \approx \dfrac{1}{6\Delta x}\delta^{(1)_3}, \delta^{(1)_3} = (-f_{j+2} + 6f_{j+1} - 3f_j - 2f_{j-1}) \\ \left(\dfrac{\partial^2 f}{\partial x^2}\right)_j \approx \dfrac{1}{\Delta x^2}\delta^{(2)_2} : \delta^{(2)_2} = (f_{j+1} - 2f_j + f_{j-1}) \\ \left(\dfrac{\partial^3 f}{\partial x^3}\right)_j \approx \dfrac{1}{\Delta x^3}\delta^{(3)_1}, \delta^{(3)_1} = (f_{j+2} - 3f_{j+1} + 3f_j - f_{j-1}) \end{cases} \quad (14)$$

It worthwhile to mention that there are other discretizations, however, we have verified they would not contribute more or fail to yield a $\tau$ so that $\tau/\beta_0^{(2)}, \tau/\beta_2^{(3)} \sim O(\Delta x)$ at critical points. Based on Eq. (14), candidates of multiplication of two discretizations which satisfy aforementioned accuracy relations are selected and their accuracy relations at critical points shown in Table 8. Besides, it is easy to find that in the absence of critical points, the leading errors of $\{\beta_k^{(2)}, \beta_k^{(3)}, \delta^{(1)_3} \cdot \delta^{(3)_1}, \delta^{(2)_2} \cdot \delta^{(3)_1}, {\delta^{(3)_1}}^2\}$ are of the quantity as $\{\Delta x^2, \Delta x^2, \Delta x^4, \Delta x^5, \Delta x^6\}$, and therefore Eq. (8) is satisfied if $\tau$ takes the corresponding multiplication.

Table 8 Accuracy relations of $\beta_k^{(3)}$, multiplications of derivative discretizations and $\tau_{CP_1}$ in occurrences of CP$_0$ and CP$_1$

| | $\lambda$ | CP$_1$ | CP$_2$ |
|---|---|---|---|
| $\beta_0^{(3)}$ | $\neq 0$ | $\left(\dfrac{13}{12}+\lambda^2\right){f''_{x_c}}^2\Delta x^4 - (\lambda^3 + \dfrac{3}{2}\lambda + \dfrac{13}{6})f''_{x_c}f'''_{x_c}\Delta x^5 + O(\Delta x^6)$ | $(\dfrac{1}{4}\lambda^4 + \dfrac{3}{4}\lambda^2 + \dfrac{13}{6}\lambda + \dfrac{43}{36}){f'''_{x_c}}^2\Delta x^6 - (\dfrac{1}{6}\lambda^5 + \dfrac{23}{36}\lambda^3 + 3\lambda^2 + \dfrac{263}{72}\lambda + \dfrac{103}{72})f'''_{x_c}f^{(4)}_{x_c}\Delta x^7 + O(\Delta x^8)$ |
| | $= 0$ | $\dfrac{13}{12}{f''_j}^2\Delta x^4 - \dfrac{13}{6}f''_jf'''_j\Delta x^5 + O(\Delta x^6)$ | $\dfrac{43}{36}{f'''_j}^2\Delta x^6 - \dfrac{103}{72}f'''_jf^{(4)}_j\Delta x^7 + O(\Delta x^8)$ |
| $\beta_2^{(3)}$ | $\neq 0$ | $\left(\dfrac{13}{12}+\lambda^2\right){f''_{x_c}}^2\Delta x^4 - (\lambda^3 + \dfrac{3}{2}\lambda - \dfrac{13}{6})f''_{x_c}f'''_{x_c}\Delta x^5 + O(\Delta x^6)$ | $(\dfrac{1}{4}\lambda^4 + \dfrac{3}{4}\lambda^2 - \dfrac{13}{6}\lambda + \dfrac{43}{36}){f'''_{x_c}}^2\Delta x^6 - (\dfrac{1}{6}\lambda^5 + \dfrac{23}{36}\lambda^3 - 3\lambda^2 + \dfrac{263}{72}\lambda - \dfrac{103}{72})f'''_{x_c}f^{(4)}_{x_c}\Delta x^7 + O(\Delta x^8)$ |
| | $= 0$ | $\dfrac{13}{12}{f''_j}^2\Delta x^4 + \dfrac{13}{6}f''_jf'''_j\Delta x^5 + O(\Delta x^6)$ | $\dfrac{43}{36}{f'''_j}^2\Delta x^6 + \dfrac{103}{72}f'''_jf^{(4)}_j\Delta x^7 + O(\Delta x^8)$ |
| $\delta^{(1)_3} \cdot \delta^{(3)_1}$ | $\neq 0$ | $-6\lambda f''_{x_c}f'''_{x_c}\Delta x^5 + 3((2\lambda^2-\lambda)f''_{x_c}f^{(4)}_{x_c} + \lambda^2 {f'''_{x_c}}^2)\Delta x^6 + O(\Delta x^7)$ | $3\lambda^2 {f'''_{x_c}}^2\Delta x^6 - (4\lambda^3 - \dfrac{3}{2}\lambda^2 + \dfrac{1}{2})f'''_{x_c}f^{(4)}_{x_c}\Delta x^7$ |
| | $= 0$ | $-\dfrac{1}{2}f'''_jf^{(4)}_j\Delta x^7 - (\dfrac{1}{4}{f^{(4)}_j}^2 +$ | $-\dfrac{1}{2}f'''_jf^{(4)}_j\Delta x^7 - (\dfrac{1}{4}{f^{(4)}_j}^2 +$ |
| | | $\frac{1}{5}f_j'''f_j^{(5)})\Delta x^8 + O(\Delta x^9)$ | $\frac{1}{5}f_j'''f_j^{(5)})\Delta x^8$ |
| $\delta^{(2)_2} \cdot \delta^{(3)_1}$ | $\neq 0$ | $f_{x_c}''f_{x_c}'''\Delta x^5 + ((-\lambda + \frac{1}{2})f_{x_c}''f_{x_c}^{(4)} - \lambda f_{x_c}'''^2)\Delta x^6 + O(\Delta x^7)$ | $-\lambda f_{x_c}'''^2\Delta x^6 - (\frac{3}{2}\lambda^2 - \frac{1}{2}\lambda + \frac{1}{12})f_{x_c}''f_{x_c}^{(4)}\Delta x^7$ |
| | $= 0$ | $f_j''f_j'''\Delta x^5 + \frac{1}{2}f_j''f_j^{(4)}\Delta x^6 + O(\Delta x^7)$ | $\frac{1}{12}f_j'''f_j^{(4)}\Delta x^7 + \frac{1}{24}f_j^{(4)2}\Delta x^8 + O(\Delta x^9)$ |
| $\delta^{(3)_1}{}^2$ | $\neq 0$ | $f_{x_c}'''^2\Delta x^6 - (2\lambda - 1)f_{x_c}'''f_{x_c}^{(4)}\Delta x^7 + O(\Delta x^8)$ | $f_{j_c}'''^2\Delta x^6 - (2\lambda - 1)f_{x_c}'''f_{x_c}^{(4)}\Delta x^7 + O(\Delta x^8)$ |
| | $= 0$ | $f_j'''^2\Delta x^6 + f_j'''f_j^{(4)}\Delta x^7 + O(\Delta x^8)$ | $f_j'''^2\Delta x^6 + f_j'''f_j^{(4)}\Delta x^7 + O(\Delta x^8)$ |
| $\tau_{CP_1}$ | $\neq 0 \& -\frac{1}{2}$ | $(-6\lambda - 3)f_{x_c}''f_{x_c}'''\Delta x^5 + ((6\lambda^2 - \frac{3}{2})f_{x_c}''f_{x_c}^{(4)} + (3\lambda^2 + \frac{3}{4})f_{x_c}'''^2)\Delta x^6 + O(\Delta x^7)$ | -- |
| | $= -\frac{1}{2}$ | $-\frac{1}{4}f_j'''f_j^{(4)}\Delta x^7 + O(\Delta x^8)$ | -- |
| | $= 0$ | $-3f_j''f_j'''\Delta x^5 + (\frac{3}{4}f_j'''^2 - \frac{3}{2}f_j''f_j^{(4)})\Delta x^6 + O(\Delta x^7)$ | -- |

According to above discussion especially the Table 6 and 9, one can see that if $\tau$ takes the one of $\{\delta^{(1)_3} \cdot \delta^{(3)_1}, \delta^{(2)_2} \cdot \delta^{(3)_1}, \delta^{(3)_1}{}^2\}$, $\tau/\beta_k^{(2)}$ would make Eq. (8) satisfied if $\lambda \neq \pm 1/2$. However, when $\lambda = \pm 1/2$, $\beta_k^{(2)} \sim O(\Delta x^6)$, and therefore even if $\beta_1^{(2)}$ is replaced by the extension $\beta_2^{(3)}$, Eq. (8) cannot be fulfilled. To solve this problem, the following measures are taken: (1) Upgrade $\beta_1^{(2)}$ with $\beta_2^{(3)}$ as previously suggested, such that the case $\lambda = 1/2$ won't cause trouble; (2) Try combining $\{\delta^{(1)_3} \cdot \delta^{(3)_1}, \delta^{(2)_2} \cdot \delta^{(3)_1}, \delta^{(3)_1}{}^2\}$ to produce a $\tau$ such that when $\lambda = 1/2$, $\tau$ would have the error as $O(\Delta x^7)$, and therefore $\tau/\beta_0^{(2)}$ would have the error as $O(\Delta x)$. Fortunately, the solution exists which takes the form as:

$$\tau_{CP_1} = \left|\delta^{(1)_3} \cdot \delta^{(3)_1} - 3\delta^{(2)_2} \cdot \delta^{(3)_1} + \frac{3}{4}\delta^{(3)_1}{}^2\right|$$
$$= c \times \left|(-f_{j+2} + 3f_{j+1} + 21f_j - 23f_{j-1}) \times (f_{j+2} - 3f_{j+1} + 3f_j - f_{j-1})\right| \quad (15)$$

where $c = 1/4$ and the subscript "CP$_1$" denotes the capability of optimal order recovery at the critical point with the first-order. In detail, the accuracy relations of $\tau_{CP_1}$ are tabulated in Table 8, from which and Table 8, one can verify the establishment of $\tau_{CP_1}/\beta_0^{(2)}, \tau_{CP_1}/\beta_2^{(3)} \sim O(\Delta x^n)$ with $n \geq 1$ at the first-order critical point. Hence, by means of formulation similar to Eq. (6), namely $\alpha_k = d_k(1 + c(\tau/\beta_k)^2)$ where $0 < c \leq 1$ and $\tau/\beta_k$ takes $\tau_{CP_1}/\beta_0^{(2)}$ and $\tau_{CP_1}/\beta_2^{(3)}$ respectively, $\alpha_k$ can theoretically make Eq. (8) satisfied not only at non-critical point situation but also at the first-order critical points. Furthermore, we propose the following proposition to indicate Eq. (15) with the inside coefficient 1/4 replaced as $c$ is only solution fulfilling the required accuracy relation at the first-order critical point.

**Proposition 6.** Consider the generic quadratic form of $f$ as $\tau(f) = (f_{j-1}, f_j, f_{j+1}, f_{j+2})[a_{i_1 i_2}](f_{j-1}, f_j, f_{j+1}, f_{j+2})^T$ where $[a_{i_1 i_2}]$ is a $4 \times 4$ matrix and $i_1, i_2 = 1\ldots 4$. Supposing the first-order critical point would occur at $x_c = x_j + \lambda \cdot \Delta x$ where $-1 < \lambda < 2$ and $\{f'_{x_c} = 0, f''_{x_c} \& f'''_{x_c} \neq 0\}$, the solution of $\tau(f)$ is wanted such that its Taylor expansion toward $x_c$ having the leading error as: (1) $O(\Delta x^5)$ when $\lambda \neq 0 \& -\frac{1}{2}$; (2) $O(\Delta x^7)$ when $\lambda = -\frac{1}{2}$, then the only solution is found as: $\tau(f) = c \times (-f_{j+2} + 3f_{j+1} + 21f_j - 23f_{j-1}) \times (f_{j+2} - 3f_{j+1} + 3f_j - f_{j-1})$ where $c$ is a free parameter.

The proof of the proposition is given in Appendix II for brevity.

From the above analysis, it can be seen that the formulation $\alpha_k = d_k(1 + c(\tau/\beta_k)^p)$ with $p = 2$ can theoretically make corresponding WENO-Z scheme achieve the third-order at the first-order critical point and be SCI. However, the following dissatisfactions exist: (1) Proposition 4 indicates $p = 2$ would make the scheme have less resolution on small structures than that with $p = 1$, as already shown in Ref [4]; (4) Proposition 3 indicates that the smaller $c$ will help the scheme to increase resolution and $\alpha_k$ be more close to $d_k$ at smooth region, however, inappropriately small $c$ would have the risk of numerical instability near discontinuities. Hence, a dynamic application of $p$ and $c$ are expected.

As referred in the introduction, WENO-M is another representative improvement of WENO-JS. Henrick et al. [3] first indicated that if one mapping function regarding $\omega_k$ and $d_k$ in Eq. (4) exists which satisfies: $g(d_k; \omega_k) = d_k + O(\omega_k - d_k)^n$ withe $n$ as certain positive integer and $n \geq 2$ usually, then as long as $\omega_k - d_k = O(\Delta x^{n_1})$ with $n_1 \geq 1$, $g(d_k; \omega_k) = d_k + O(\Delta x^{n \times n_1})$ holds where $g(d_k; \omega_k)$ acts as the new $\alpha_k$. One can see Eq. (8) would be satisfied thereafter and the optimal order of WENO-JS would be recovered in the occurrence of critical points. In Ref. [4], a concrete mapping was proposed for WENO5-JS, and subsequent developments were also observed. Recently, we made intensive studies [14] on WENO-M methodology and proposed a piecewise rational mapping which has fine regulation capability. In order to facilitate the discussion on the mapping characteristics, we propose the so-called $C_{n,m,l}$ condition in the region [0, 1], and the piecewise version $C_{n,m}$ in $[d_k, 1]$ can be simplified as:

$$g^{(i)}(d_k) = \begin{cases} d_k, & i = 0 \\ 0, & 1 \leq i \leq n \\ \neq 0, & i = n+1 \end{cases}, \quad g^{(i)}(1) = \begin{cases} 1, & i = 0 \\ 1, & i = 1 \text{ if } l \geq 1 \\ 0, & 2 \leq i \leq l \\ \neq 0, & i = l+1 \end{cases}$$

The proposed piecewise rational mapping in $[d_k, 1]$ which satisfies $C_{n,\min(m,m_1-1)}$ is [14]

$$PRM_{n,m;m_1;c_1,c_2}^{n+1} = \frac{(\omega - d_k)^{n+1}}{(\omega - d_k)^n + c_2(\omega - d_k)(1-\omega)^{m_1} + c_1(1-\omega)^{m+1}}$$

In order to apply the above outcomes for aforementioned dissatisfactions, we take the same thinking and extend the above $C_{n,m}$ as

$$M^{(i)}(0) = \begin{cases} 0, & 1 \leq i \leq n \\ \neq 0, & i = n+1 \end{cases}, \quad M(c_3) = \begin{cases} c_3, & i = 0 \\ 1, & i = 1 \text{ if } l \geq 1 \\ 0, & 2 \leq i \leq l \\ \neq 0, & i = l+1 \end{cases} \quad (15)$$

where $M(\cdot)$ is the mapping faction and $c_3$ is one positive parameter prescribed later. Then, as along with Eq. (15), original $PRM_{n,m;m_1;c_1,c_2}^{n+1}$ is extended as

$$M_{n,m;m_1;c_1,c_2,c_3}^{n+1}(\varpi) = \begin{cases} \dfrac{\varpi^{n+1}}{\varpi^n + c_2\varpi(c_3-\varpi)^{m_1} + c_1(c_3-\varpi)^{m+1}}, & \varpi \leq c_3 \\ \varpi, & \varpi > c_3 \end{cases}, \quad (16)$$

Where $m_1 \geq m+1$. Considering the satisfaction of $C_{n,\min(m,m_1-1)}$ by $PRM_{n,m;m_1;c_1,c_2}^{n+1}$ [11], it is easy to see that Eq. (16) would make Eq. (15) hold.

Combining Eq. (16) with $\alpha_k = d_k(1 + c(\tau/\beta_k)^p)$, a new non-normalized weight is defined as:

$$\alpha_k = d_k(1 + M(\tau/\beta_k)). \quad (17)$$

where the mapping $M(\cdot)$ would take Eq. (16) is this study. From Eq. (15), one can see that as long as $n > 2$ and $\tau/\beta_k \sim \Delta x^{n_1}$ with $n_1 \geq 1$, $M(\tau/\beta_k) \sim O(\Delta x^{>2n_1})$ holds and Eq. (8) is established, while for $\tau/\beta_k \sim O(1)$ and above, which represents the occurrence of discontinuities, $\alpha_k \sim d_k(1 + (\tau/\beta_k))$ which resembles the performance of WENO5-Z.

Specifically, for current WENO3 scheme, the mapping is chosen as $M_{2,1}^3$ where $m_1 = m+1$ and the parameters $\{c_1, c_2, c_3\}$ are defined with respect to $d_k$ after extensive numerical practices. The results are shown in Table 9 and the mapping distributions are shown in Fig. 3.

Table 9 Coefficients of $M_{2,1;2;c_1,c_2,c_3}^3$ in Eq. (16) with respect to $d_k$

|  | $c_1$ | $c_2$ | $c_3$ |
|---|---|---|---|
| $d_0 = 1/3$ | 1.2 | 0.1 | 55 |
| $d_1 = 2/3$ | 1.2 | 0.1 | 35 |

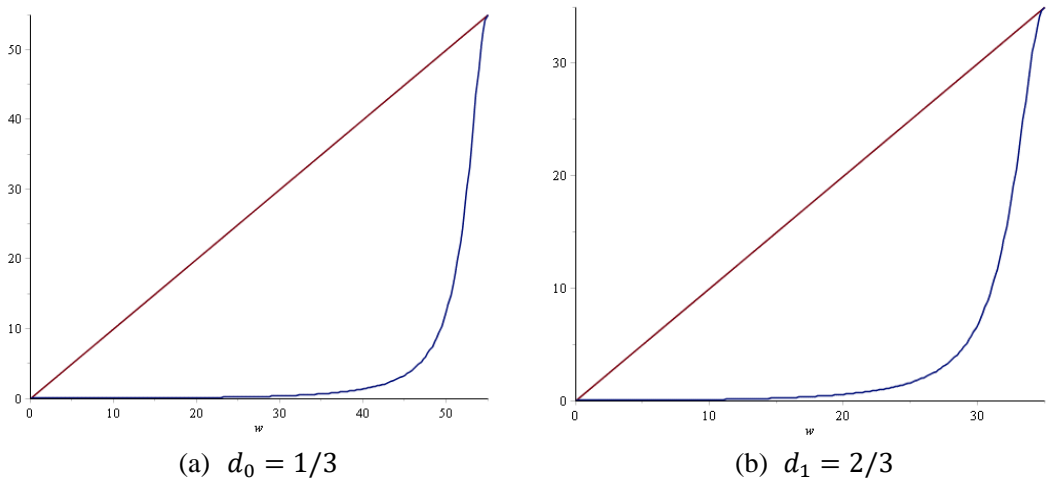

(a) $d_0 = 1/3$          (b) $d_1 = 2/3$

Fig. 3 The distributions of $M_{2,1}^3$ in cases of $d_0 = 1/3$ and $d_1 = 2/3$

So far, a third-order WENO-Z improvement is accomplished which can recover the optimal order at $\{f'_{x_c} = 0, f''_{x_c} \neq 0, f'''_{x_c} \neq 0\}$ and is SCI. For reference, the scheme is referred as WENO3-ZM where "M" denotes mapping; to facilitate programming, its implementation is summarized as:

(1) Compute $\beta_k$ by: $\beta_0 = \beta_0^{(2)}$ and $\beta_1 = \beta_2^{(3)}$, where $\beta_k^{(r)}$ is defined by Eq. (5).
(2) Compute $\tau$ by: $\tau = \tau_{cp1}$ from Eq. (15).
(3) Compute $\alpha_k$ by Eq. (17) where $M(\cdot)$ is defined by Eq. (16) and corresponding parameters are given in Table 9.
(4) Eqns. (3) -(4) are implemented to acquire the final $\hat{f}_{j+1/2}$.

4.2 New third-order WENO-Z scheme by further expanding stencil

In the section just above, the practice is made by only extending $\beta_1^{(2)}$ to $\beta_2^{(3)}$, where the global stencil is expanded accordingly and the new $\tau_{CP_1}$ is proposed. One may naturally wonder the attempt to further extend $\beta_0^{(2)}$ to $\beta_0^{(3)}$, as that has been done in Ref. [14]. For the practice, the possible advantages include: (1) As shown in Table 8, the coefficient of the leading error term of $\beta_0^{(3)}$ will be non-zero regardless of $\lambda$, therefore the extensions totally remove the dilemma caused by the occurrence of first-order critical points at $\lambda = \pm 1/2$; (2) The ensued expansion of global stencil implies the availability of new $\tau$ which would have Eq. (8) established. In order to devise such a $\tau$, we follow the similar idea in Section 4.1 and first check the available derivatives at $x_j$ and their discretizations on $\{x_{j-2}, x_{j-1}, x_j, x_{j+1}, x_{j+2}\}$ with the possible highest order as follows.

$$\begin{cases} \left(\frac{\partial f}{\partial x}\right)_j \approx \frac{1}{12\Delta x}\delta^{(1)_4}, \delta^{(1)_4} = (-f_{j+2} + 8f_{j+1} - 8f_{j-1} + f_{j-2}) \\ \left(\frac{\partial^2 f}{\partial x^2}\right)_j \approx \frac{1}{12\Delta x^2}\delta^{(2)_4} : \delta^{(2)_4} = (-f_{j+2} + 16f_{j+1} - 30f_j + 16f_{j-1} - f_{j-2}) \\ \left(\frac{\partial^3 f}{\partial x^3}\right)_j \approx \frac{1}{2\Delta x^3}\delta^{(3)_2}, \delta^{(3)_2} = (f_{j+2} - 2f_{j+1} + 2f_{j-1} - f_{j-2}) \\ \left(\frac{\partial^4 f}{\partial x^4}\right)_j \approx \frac{1}{\Delta x^2}\delta^{(4)_2} : \delta^{(4)_2} = (f_{j+2} - 4f_{j+1} + 6f_j - 4f_{j-1} + f_{j-2}) \end{cases} \quad (18)$$

As mentioned before, there are other discretizations regarding the above derivatives, and we find they usually cannot provide more help on devising a $\tau$ which would have high-order accuracy in occurrences of critical points. Hence, corresponding practices are not described for brevity. In terms of the dimension requirement, the one of $[\tau]$ should be $[f]^2$, and therefore the multiplications of two $\delta^{(m)_n}$ in Eq. (18) are worthy of consideration. Because the expansion at this time has obviously increased computational cost, we only consider the scheme which would employ Eq. (6) with $p = 1$.

Considering the accuracy relations of $\beta_0^{(3)}$ and $\beta_2^{(3)}$ in Table 8, in order to recover the third-order at the first-order critical point, the error of $\tau$ should have the quantity as $O(\Delta x^{\geq 6})$; to fulfill the job at the second-order critical point, the error of $\tau$ should be as $O(\Delta x^{\geq 8})$. One can check that in the occurrence of the first-order critical point, $\{\delta^{(3)_2}{}^2, \delta^{(3)_2} \cdot \delta^{(4)_2}, \delta^{(4)_2}{}^2\}$ would satisfy the requirement; in the occurrence of second-order critical point, only $\delta^{(4)_2}{}^2$ satisfies the requirement. For illustration, the accuracy relations of $\delta^{(4)_2}{}^2$ are tabulated in Table 9.

Table 9 Accuracy relation of $\delta^{(4)_2}{}^2$ in occurrences of CP$_0$, CP$_1$, CP$_2$

| Critical point | Accuracy relation |
| --- | --- |
| CP$_0$ | $f_j^{(4)2}\Delta x^8 - \frac{1}{3}f_j^{(4)}f_j^{(6)}\Delta x^{10} + O(\Delta x^{12})$ |

| CP$_1$, CP$_2$ | $\lambda \neq 0$ | $f_{x_c}^{(4)2}\Delta x^8 - 2\lambda f_{x_c}^{(4)}f_{x_c}^{(5)}\Delta x^9 + O(\Delta x^{10})$ |
| --- | --- | --- |
| | $\lambda = 0$ | $f_j^{(4)2}\Delta x^8 + \frac{1}{3}f_j^{(4)}f_j^{(6)}\Delta x^{10} + O(\Delta x^{12})$ |

Thus far, a new global indicator is acquired as

$$\tau_{CP_2} = c \times \left(f_{j+2} - 4f_{j+1} + 6f_j - 4f_{j-1} + f_{j-2}\right)^2 \qquad (19)$$

where $c$ is a positive parameter and the subscript "CP$_2$" denotes the capability of optimal order recovery at the critical point with orders up to 2. When the non-normalized weight of WENO3 is defined as: $\left\{\alpha_0 = d_0\left(1 + \tau/\beta_0^{(3)}\right), \alpha_1 = d_1\left(1 + \tau/\beta_2^{(3)}\right)\right\}$, corresponding WENO-Z scheme can achieve the third-order at situations where non-critical points, first- and second-order critical points exist. According to accuracy relations in Table 8-9, one can see that the first-order critical point is referred in the sense of $\{f'_{x_c} = 0, f''_{x_c} \neq 0\}$, and the second-order one is referred at $\{f'_{x_c} = f''_{x_c} = 0, f'''_{x_c} \neq 0\}$. For reference, we refer the scheme as WENO3-Z$_{ES}$, where the subscript "ES" represents the manipulation of expending stencil, and in addition, WENO3-Z$_{ES}$ is SCI. One can see that $\alpha_k$ would have much small error with respect to $d_k$ at CP$_0$ and CP$_1$, which implies the high accuracy of WENO-Z$_{ES}$ therein; in the meanwhile, the risk of less robustness arises as shown in Section 5.

## 5 Numerical examples

5.1 Case descriptions

Three kinds of equations and corresponding problems are provided to illustrate the properties and performances of WENO3-ZM and WENO3-Z$_{ES}$, i.e. 1-D scalar advection equation, 1-D Euler equations, 2-D Euler equations.

(1) 1-D scalar advection equation

The governing equation is: $\partial u/\partial t + \partial u/\partial x = 0$ with various initial conditions $u(x, 0)$ corresponding to specific problems. Problems include:

(a) Sinusoidal-like wave advection

$$u(x, 0) = \sin\left(\pi(x - x_c) - \frac{\sin(\pi(x-x_c))}{\pi}\right) \qquad (20)$$

where $x \in [-1, 1]$ and $x_c = 0.5966831869112089637212$

In this case, $u(x, 0)$ have two first-order critical points which locates at $x = 0$ and $x = -2 + 2x_c$, and $x = \pm x_c$ are the locations of critical points of $\sin(\pi x - \sin(\pi x)/\pi)$. Since the fourth-order Runge-Kutta (RK4) scheme is used for time discretization, the requirement of time step should be: $\Delta t < \Delta x^{\frac{3}{4}}$. Considering $\Delta x < \Delta x^{\frac{3}{4}}$ when $\Delta x < 1$, we can use CFL number to define $\Delta t$ as $\Delta t = CFL \cdot \Delta x$. A series of grid cells with the numbers {10, 20, 40, 80 …} are employed and $x = 0$ coincides with certain grid point initially. The computation runs for $t = 2$, and two CFL numbers are chosen as 0.25 and 0.4. For the first choice, one of critical points will theoretically move to half node ($\lambda = 1/2$) every four times of computation step; for the second choice, the critical points will move within grid cells but do not locate on half nodes. It is worth mentioning the derivatives of critical points have the properties as: $f' = 0$, $f'' \neq 0$, and $f''' \neq 0$.

(b) Combination-waves advection

$$u(x,0) = \begin{cases} \frac{1}{6}(G(x,\beta,z-\delta) + G(x,\beta,z+\delta) + 4G(x,\beta,z)), & -0.8 \leq x \leq -0.6 \\ 1, & -0.4 \leq x \leq -0.2 \\ 1 - |10(x-0.1)|, & 0 \leq x \leq 0.2 \\ \frac{1}{6}(F(x,\alpha,a-\delta) + F(x,\alpha,a+\delta) + 4F(x,\alpha,a)), & 0.4 \leq x \leq 0.6 \\ 0, & otherwise \end{cases} \quad (21)$$

where $x \in [-1,1]$, $G(x,\beta,z) = e^{-\beta(x-z)^2}$, $F(x,\alpha,a) = \sqrt{\max(1-\alpha^2(x-a)^2, 0)}$, $a = 0.5$, $z = -0.7$, $\delta = 0.005$, $\alpha = 10$ and $\beta = \log 2/36\delta^2$. The grid number is $N = 800$, and computation advances to a long period $T = 4000$. In computations, the third-order TVD Runge-Kutta method [2] (TVD-RK3) is employed with the CFL number as 0.1.

Although this example is a canonical test which often runs in a short time period, the stability in the long-time computation is seldom concerned by investigations. As shown in Ref. [14], unexpected outcomes usually arise among various schemes. According to our knowledge, the schemes in Table 1 have not experienced such test thoroughly and corresponding outcomes are unclear.

(2) 1-D Euler equations

Because the routine computations such as Sod problems are easily accomplished by schemes, the following three problems are chosen, i.e. strong shock wave, blast wave and Shu-Osher problem.

(a) Strong shock wave

This case poses a test regarding the computational robustness. The initial condition is: $(\rho, u, p) = \begin{cases} (1,0,0.1PR), & -5 \leq x < 0 \\ (1,0,0.1), & 0 < x \leq 5 \end{cases}$ with $PR = 10^6$ currently. The computation advances to $t=0.01$, and the grid number is $N=200$.

(b) Blast wave

This canonical example serves as another trial regarding the computational robustness. The initial condition is: $(\rho, u, p) = \begin{cases} (1,0,1000), & 0 \leq x < 0.1 \\ (1,0,0.01), & 0.1 \leq x \leq 0.9 \\ (1,0,100), & 0.9 < x \leq 1 \end{cases}$, and solid wall condition is casted at boundaries. The grid number is $N=200$ and the computation advances to $t=0.038$. By convention, a result on 15001 grids by WENO5 is regarded as the "Exact" solution for reference.

(c) Shu-Osher problem

This problem is a benchmark test on the numerical resolution. The initial condition is: $(\rho, u, p) = \begin{cases} (3.857143, 2.629369, 10.3333), & -5 \leq x < -4 \\ (1 + 0.2\sin(5x), 0, 1), & -4 < x \leq 5 \end{cases}$. The computation advances to $t=1.8$. By convention, the result of WENO5 at 10001 grids is regarded as the "Exact" solution for reference.

In 1-D Euler equations, the temporal scheme applies TVD-RK3, and flux splitting uses the Steger-Warming scheme. To mitigate numerical oscillations, the characteristic variables are used in schemes.

(3) 2-D Euler equations

The following typical problems are chosen which regards resolution and robustness: 2-D Riemann problem and double Mach reflection.

(a) 2-D Riemann problem

The problem is solved in a domain: $[0,1] \times [0,1]$ which is partitioned into four parts by lines $x = 0.8$ and $y = 0.8$. The initial conditions are:

$$(\rho, u, v, p) = \begin{cases} (1.5, 0, 0, 1.5), & 0.8 \leq x \leq 1, 0.8 \leq y \leq 1 \\ (0.5323, 1.206, 0, 0.3), & 0 \leq x < 0.8, 0.8 \leq y \leq 1 \\ (0.138, 1.206, 1.206, 0.029), & 0.8 \leq x < 0.8, 0 \leq y < 0.8 \\ (0.5323, 0, 1.206, 0.3), & 0.8 \leq x \leq 1, 0 \leq y < 0.8 \end{cases}$$

The grid number is $960 \times 960$. The computation advances to $t = 0.8$ with the specific heat ratio as $\gamma = 1.4$

(b) Double Mach reflection

The problem describes a Mach 10 shock impinging the wall at an incident angle of 60°. The computational domain is $[0,3] \times [0,1]$, and initial condition is

$$(\rho, u, v, p) = \begin{cases} (8, 7.145, -4.125, 116.5), & x < 1/6 + y/\sqrt{3} \\ (1.4, 0, 0, 1), & x \geq 1/6 + y/\sqrt{3} \end{cases}$$

The flow takes $\gamma$ as 1.4. The computation runs on $1920 \times 480$ grids until $t = 0.2$.

In 2-D Euler equations, the temporal scheme applies TVD-RK3 due to the unsteady nature of problems, while the flux splitting employs the Steger-Warming scheme and the characteristic variables are used in schemes.

5.2 1-D scalar advection equation

The computations here have two purposes: the first is to illustrate the order convergence property of WENO3-ZM and -$Z_{ES}$, the second is to show their stability of long-time computation. Besides computations of WENO3-ZM and -$Z_{ES}$, that of some schemes in Table 1 are also carried out for comparison.

(1) Sinusoidal-like wave advection by Eq. (20)

We first check the numerical error and corresponding orders of proposed schemes at $CFL = 0.4$ in Table 10, which corresponds to first-order critical points occurring within intervals but not on grid nodes. Regarding the choices of comparatives, because we have manifested WENO-F3 recovers the optimal order in Table 5, and the scheme in indicated in Section 2 having small difference with WENO-NP3, the results of the two schemes are omitted here; because the accomplishment of order recovery by WENO-PZ3 has been reported in Ref. [11] which is also verified by us, only the results of WENO-NN3 are provide here. For WENO-P+3, Ref. [12] has stated its failure of optimal order recovery which is validated by us also, correspond results are omitted also. Table 10 tells that the schemes can achieve the optimal third-order, which is consistent with analyses in Sections 2 and 4.

Table 10 Two norms of error and corresponding orders of WENO3-ZM, -$Z_{ES}$ and WENO-NN3 by 1-D scalar advection equation with initial condition Eq. (20) at $t = 2$ and $CFL = 0.4$

| N | Δt | WENO3-ZM | | WENO3-$Z_{ES}$ | | WENO-NN3 | |
|---|---|---|---|---|---|---|---|
| | | $L_1$-error | $L_1$-order | $L_1$-error | $L_1$-order | $L_1$-error | $L_1$-order |
| 10 | 0.08 | 9.1702E-02 | | 8.9549E-02 | | 1.2368E-01 | |
| 20 | 0.04 | 2.2331E-02 | 2.037 | 2.1600E-02 | 2.051 | 2.3591E-02 | 2.390 |
| 40 | 0.02 | 3.1837E-03 | 2.810 | 3.1621E-03 | 2.772 | 3.2592E-03 | 2.855 |
| 80 | 0.01 | 4.0473E-04 | 2.975 | 4.0451E-04 | 2.966 | 4.1244E-04 | 2.982 |
| 160 | 0.005 | 5.0827E-05 | 2.993 | 5.0819E-05 | 2.992 | 5.0858E-05 | 3.019 |
| 320 | 0.0025 | 6.3641E-06 | 2.9972 | 6.3642E-06 | 2.997 | 6.3640E-06 | 2.998 |

| 640 | 0.00125 | 7.9609E-07 | 2.998 | 7.9609E-07 | 2.998 | 7.9610E-07 | 2.998 |
|---|---|---|---|---|---|---|---|
|  |  | $L_\infty$-error | $L_\infty$-order | $L_\infty$-error | $L_\infty$-order | $L_\infty$-error | $L_\infty$-order |
| 10 | 0.08 | 2.2603E-01 |  | 2.2763E-01 |  | 2.5786E-01 |  |
| 20 | 0.04 | 5.2328E-02 | 2.110 | 5.1481E-02 | 2.144 | 5.4833E-02 | 2.233 |
| 40 | 0.02 | 7.8839E-03 | 2.730 | 7.8932E-03 | 2.705 | 6.9930E-03 | 2.971 |
| 80 | 0.01 | 1.0214E-03 | 2.948 | 1.0191E-03 | 2.953 | 1.0165E-03 | 2.782 |
| 160 | 0.005 | 1.2816E-04 | 2.994 | 1.2815E-04 | 2.991 | 1.2816E-04 | 2.987 |
| 320 | 0.0025 | 1.6035E-05 | 2.998 | 1.6035E-05 | 2.998 | 1.6034E-05 | 2.998 |
| 640 | 0.00125 | 2.0047E-06 | 2.999 | 2.0047E-06 | 2.999 | 2.0047E-06 | 2.999 |

Next, we check the situations at $CFL = 0.25$ in Table 11, which corresponds to the critical point occurring on half node occasionally. Due to the similar reason as above, only the results of WENO-NN3 are chosen as the comparison here. The table tells that all $L_1$-orders indicate the recovery of optimal order, however the failure of $L_\infty$-order by WENO-NN3 occurs. The results on the one hand manifests the true achievement of optimal order at first-order critical points by {WENO3-ZM, –$Z_{ES}$} and therefore the correctness of ACP$_\lambda$, one the other hand indicate the possible illusion in the absence of $L_1$-order analysis.

From the above discussion, one can see that all third-order improvements in Table 1 fail to recover the optimal order at first-order critical points occurring on the half node.

Table 11 Two norms of error and corresponding orders of WENO3-ZM, -$Z_{ES}$ and WENO-NN3 by 1-D scalar advection equation with initial condition Eq. (20) at $t = 2$ and $CFL = 0.25$

| N | Δt | WENO3-ZM | | WENO3-$Z_{ES}$ | | WENO-NN3 | |
|---|---|---|---|---|---|---|---|
|  |  | $L_1$-error | $L_1$-order | $L_1$-error | $L_1$-order | $L_1$-error | $L_1$-order |
| 10 | 0.08 | 9.3127E-02 |  | 8.9473E-02 |  | 1.2496E-01 |  |
| 20 | 0.04 | 2.2559E-02 | 2.045 | 2.1590E-02 | 2.051 | 2.3649E-02 | 2.401 |
| 40 | 0.02 | 3.1665E-03 | 2.832 | 3.1615E-03 | 2.771 | 3.2976E-03 | 2.842 |
| 80 | 0.01 | 4.0586E-04 | 2.963 | 4.0449E-04 | 2.966 | 4.1543E-04 | 2.988 |
| 160 | 0.005 | 5.0806E-05 | 2.997 | 5.0818E-05 | 2.992 | 5.1025E-05 | 3.025 |
| 320 | 0.0025 | 6.3639E-06 | 2.997 | 6.3641E-06 | 2.997 | 6.3630E-06 | 3.003 |
| 640 | 0.00125 | 7.9682E-07 | 2.997 | 7.9609E-07 | 2.998 | 7.9608E-07 | 2.998 |
|  |  | $L_\infty$-error | $L_\infty$-order | $L_\infty$-error | $L_\infty$-order | $L_\infty$-error | $L_\infty$-order |
| 10 | 0.08 | 2.1947E-01 |  | 2.2759E-01 |  | 2.6175E-01 |  |
| 20 | 0.04 | 5.2384E-02 | 2.066 | 5.1451E-02 | 2.145 | 5.4024E-02 | 2.276 |
| 40 | 0.02 | 7.8792E-03 | 2.733 | 7.8916E-03 | 2.704 | 6.8788E-03 | 2.973 |
| 80 | 0.01 | 1.0197E-03 | 2.949 | 1.0190E-03 | 2.953 | 1.0242E-03 | 2.747 |
| 160 | 0.005 | 1.2836E-04 | 2.989 | 1.2815E-04 | 2.991 | 1.2815E-04 | 2.998 |
| 320 | 0.0025 | 1.6035E-05 | 3.000 | 1.6035E-05 | 2.998 | 1.8132E-05 | 2.821 |
| 640 | 0.00125 | 2.0047E-06 | 2.999 | 2.0047E-06 | 2.999 | 2.8397E-06 | 2.674 |

(2) Combination-waves advection by Eq. (21)

Although this problem had been tested by WENO-NN3, -PZ3 and –P+3 [10-12], the computation only ran in short period as $t = 8$. As once reported, the issue of instability by long-time computations existed, e.g. the case was intensively investigated in Ref. [14] at $t \geq 4000$. Hence the case is studied herein with the conditions described in Section 5.1. In the study, besides the proposed schemes, the improvements in Table 1 are all tested for comparison, and the results are shown in Fig. 4.

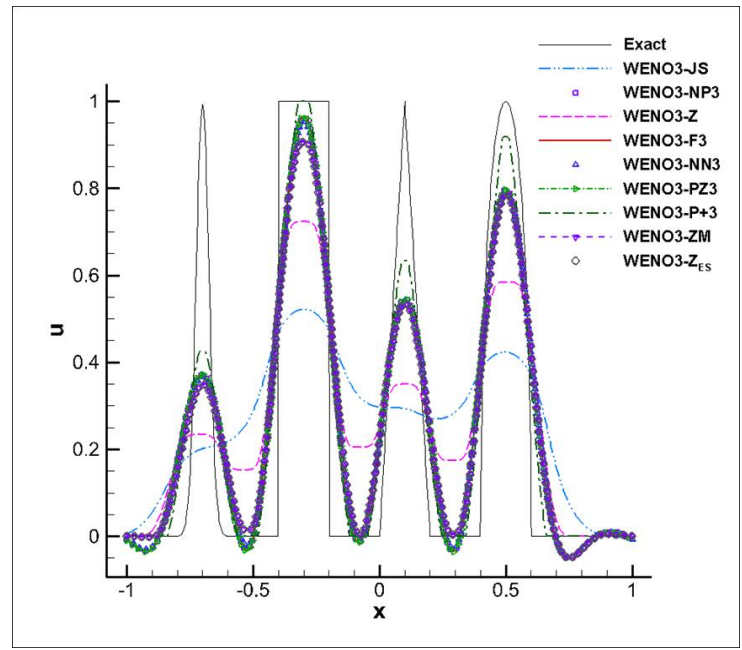

Fig. 4. The distributions of Combination-waves advection by Eq. (21) of proposed schemes with the comparisons by WENO3-JS, WENO-NP3, -F3, -NN3, -PZ3 and -P+3 at $t = 4000$

The figure tells that all improvements show oscillations except WENO-P+3, therefore the necessity of investigation cannot be neglected. Among results of the formers, the proposed WENO3-ZM and $-Z_{ES}$ show the relatively less oscillations, namely the only occurrence at the right foot of the fourth square. Although the absence of oscillations by that of WENO-P+3, its failure on cases such 1-D strong shock wave and strong oscillations in 2-D Riemann problem indicates its robustness is not fully guaranteed.

5.3 1-D Euler equations

By convention, the results of WENO3-JS are chosen as one of references automatically.

(1) Strong shock wave

In order to investigate the robustness of proposed WENO3-ZM and $-Z_{ES}$, the tests are carried out with the comparisons by WENO3-Z improvements in Table 1. It is interesting to note that WENO-P+3 fails to accomplish the computation, which indicates the relatively weak robustness. In the enlarged windows of the figure, the proposed schemes indicate comparatively better agreement with the exact solution.

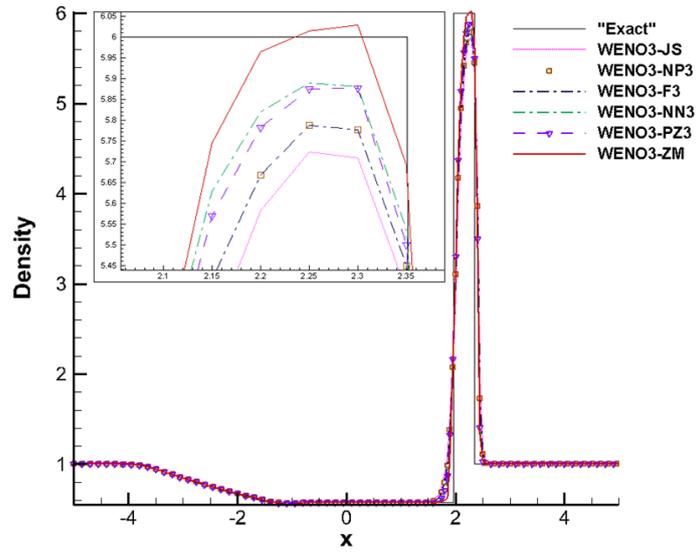

Fig. 5. Density distributions of strong shock wave at $t=0.01$ on 200 grids with initial pressure ratio $PR=10^6$ by WENO3-ZM and $-Z_{ES}$ with the comparisons by WENO3-JS, WENO-NP3, -F3, -NN3, -PZ3 and -P+3.

(2) Blast wave

The same schemes are tested as above and the results are shown in Fig. 6. Among the schemes, WENO-P+3 indicates a distribution closest to the "Exact" solution, however its failure in the previous test reminds the possiblity of inadequate numerical dissipation. Among the rest schemes, WENO3-ZM and $-Z_{ES}$ seem to show a performance with well resolution.

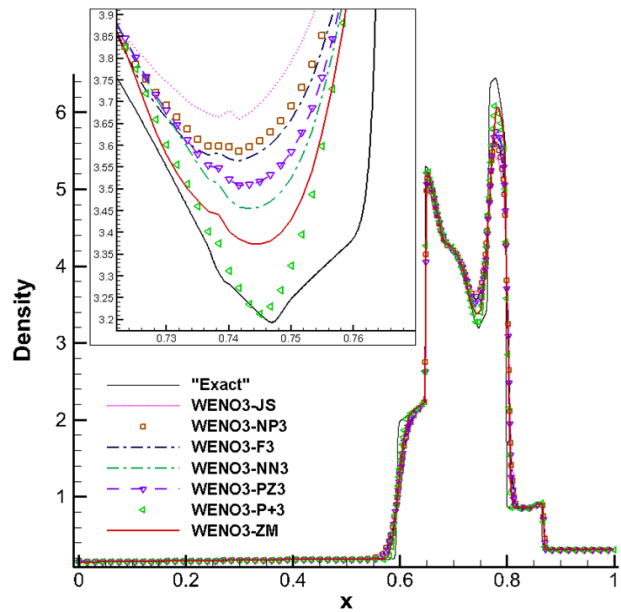

Fig. 6. Density distributions of blast waves at $t=0.038$ on 200 grids by WENO3-ZM and $-Z_{ES}$ with the comparisons by WENO3-JS, WENO-NP3, -F3, -NN3, -PZ3 and -P+3.

(3) Shu-Osher problem

Other than grids with numbers usually from 400-600 for third-order schemes, 240 grids are employed here, which poses a tough test on numerical resolution. Corresponding results are shown in Fig. 7, where $\Delta t = 0.003$ is employed in computation. From the figure, WENO3-ZM and $-Z_{ES}$ indicates resolution obviously surpassing that of WENO-NP3, -NN3 and $-PZ3$, and comparatively outperforms WENO-P+3 on the whole.

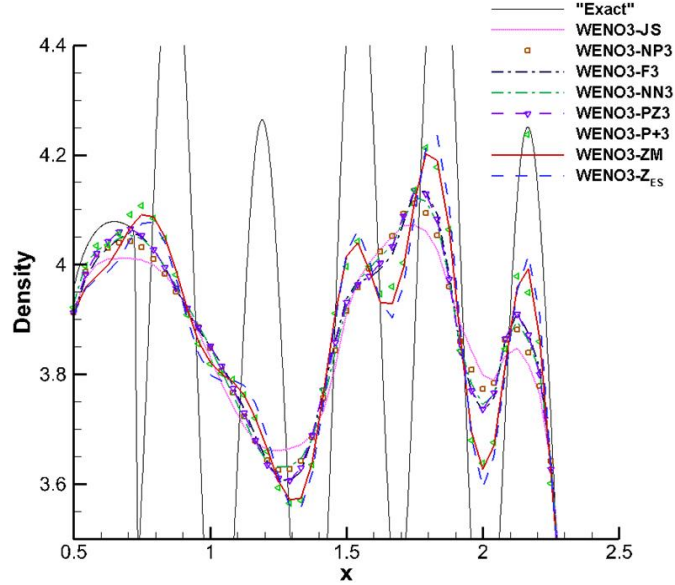

Fig. 7. Density distributions of Shu-Osher problem at $t=1.8$ & $\Delta t = 0.003$ on 240 grids by WENO3-ZM and $-Z_{ES}$ with the comparisons by WENO3-JS, WENO-NP3, -F3, -NN3, -PZ3 and -P+3.

5.4 2-D Euler equations

In 2-D problems, WENO3-JS and WENO3-Z are chosen as the basic reference scheme, and improvements {WENO-NP3, -F3, -NN3, -PZ3 and $-P+3$} are employed to compare with the proposed WENO3-ZM and $-Z_{ES}$.

(1) 2-D Riemann problem

The density contours of the schemes are illustrated in Fig. 8. Because the results of WENO-P+3 are too complicated, therefore suitable comments are not available yet. Taking the result of WENO3-JS as illustration, the following structure features in three regions 1-3 are concerns: (1) The typical Kelvin-Helmholtz instability along slip line in Region 1; (2) The evolution of the structures after the slip line interacting with the reflection shock in Region2; (3) The symmetry of structures in the region 3. Regarding the concerns: (1) For "(1)", one can see that WENO-F3, WENO3-ZM as well as WENO3-$Z_{ES}$ has produced obvious wavelets along the slip line, where the resolution indicated by WENO3-ZM seems to be relatively clearer; (2) For "(2)", the results of WENO3-ZM and WENO3-F3 indicate more resolutions; (3) For "(3)", structures in Region 3 by WENO-F3, -PZ3 indicates asymmetry to some extent, and although WENO3-JS, -Z, -NP3 engender symmetric solutions somehow, their resolutions are not striking comparatively. It seems that WENO3-ZM on the whole yields a result not only symmetric but also of well resolution.

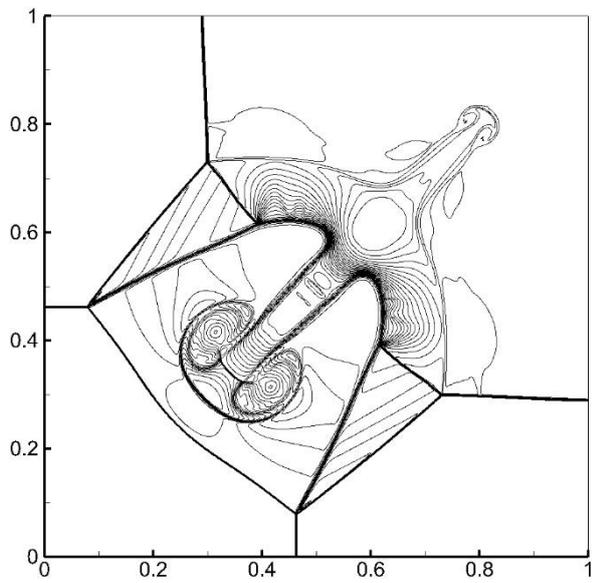
(a) WENO3-JS

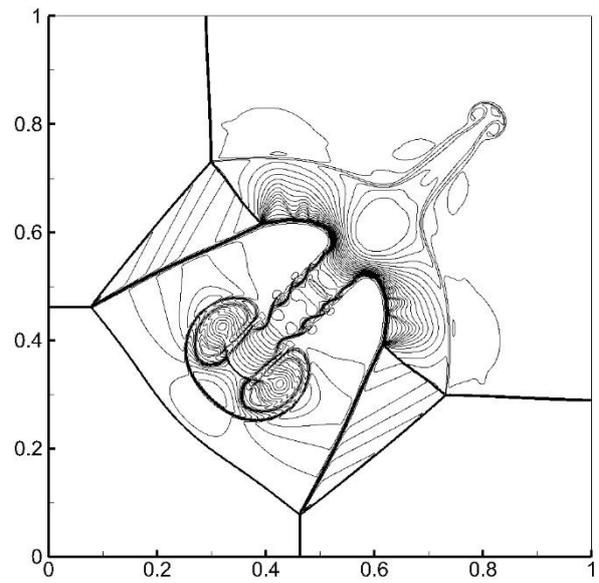
(b) WENO3-Z

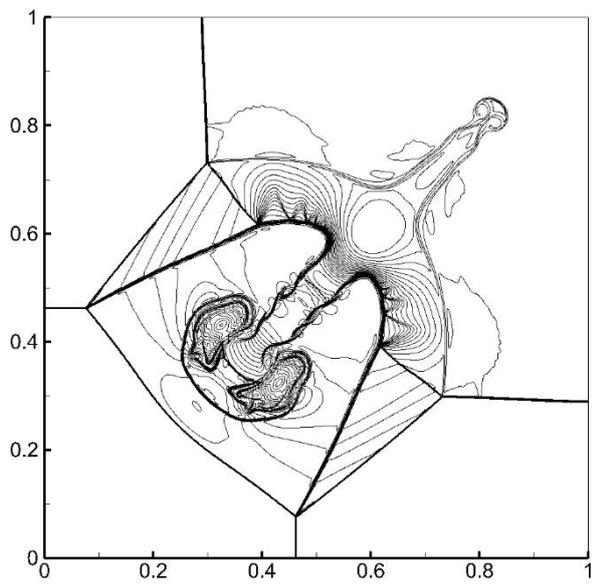
(c) WENO-NP3

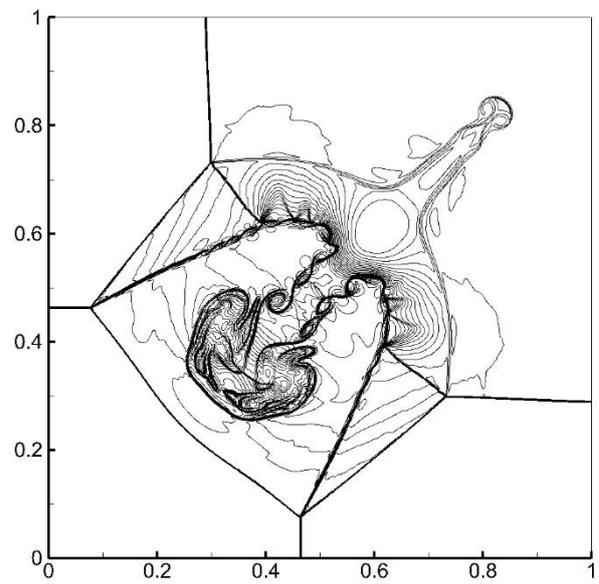
(d) WENO-F3

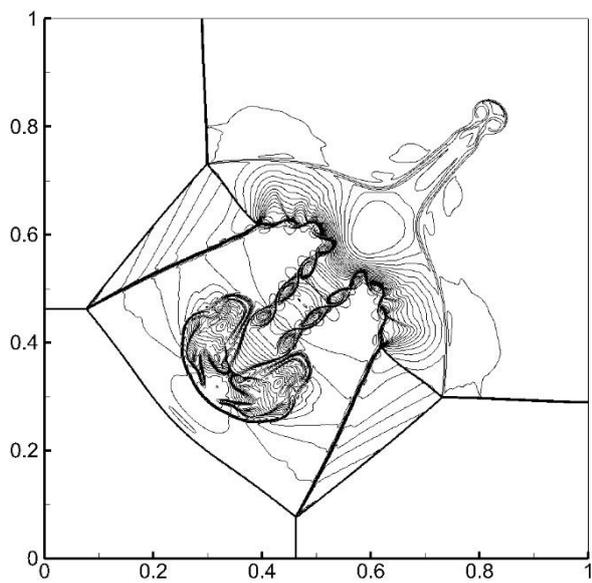

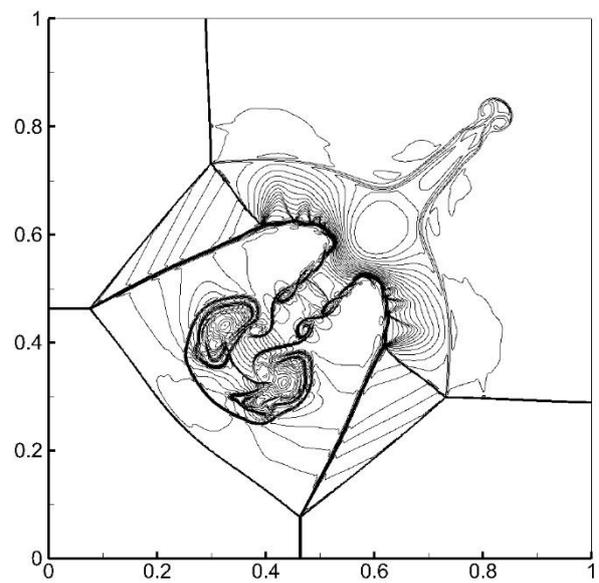

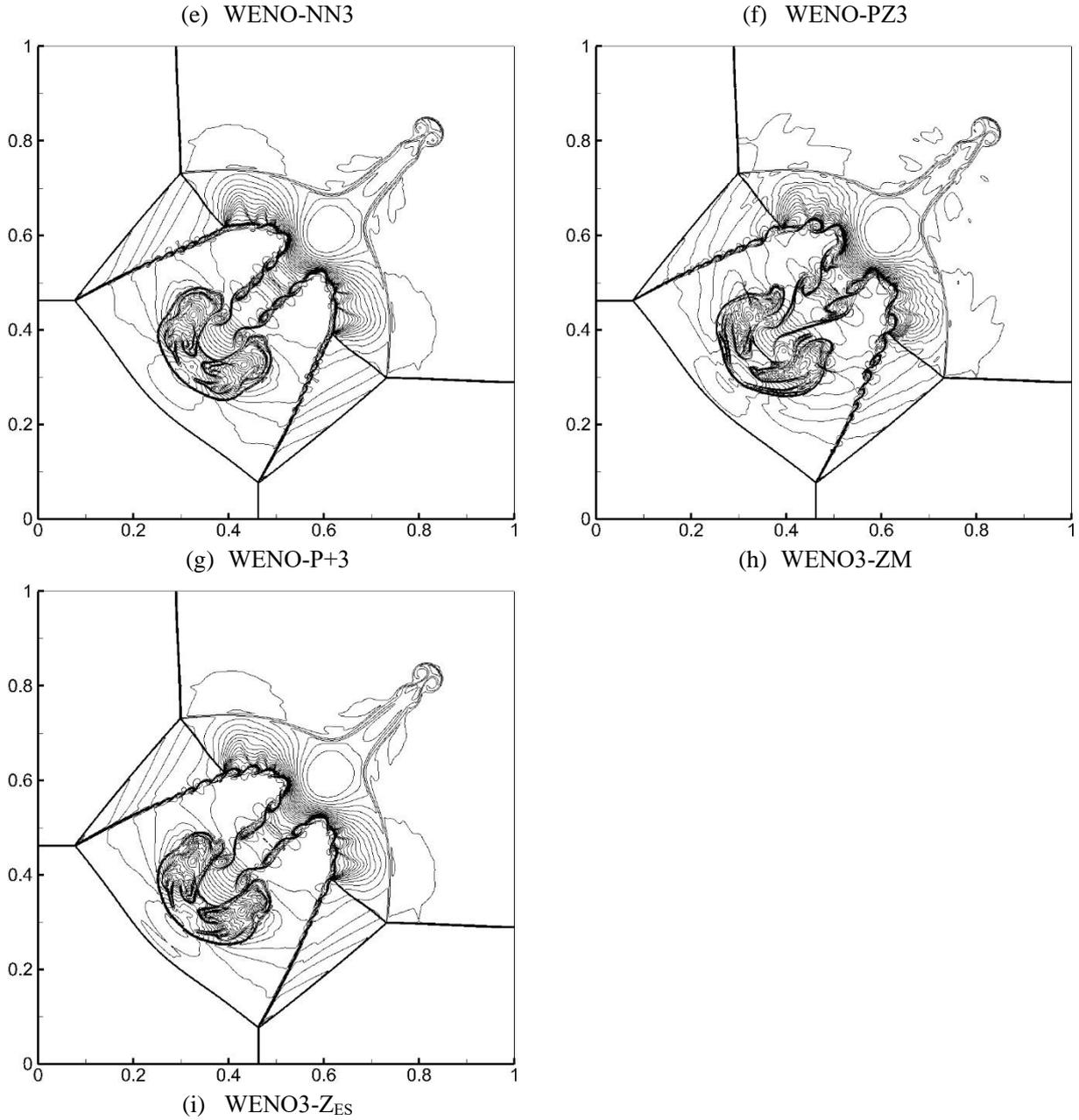

Fig. 8. Density contours of 2-D Riemann problem at on grids $240 \times 960$ at $t = 0.8$ by WENO3-ZM and –$Z_{ES}$ with the comparisons of WENO3-JS, WENO3-Z and its other improvements (40 contours from the range of 0.14 to 1.7)

(2) Double Mach reflection

It is well known that this problem tests the robustness and resolution of schemes. WENO-$Z_{ES}$ fails to accomplish the computation due to its insufficient dissipation and robustness. Corresponding results are shown in Fig. 9. It is worth mentioning we have qualitatively checked that, current results of schemes except that of WENO-F3 and WENO3-ZM are consistent with that shown in Ref. [11].

From the figure, one can see that WENO-NP3, -F3, -NN3, -PZ3 and –P+3 have shown improved resolutions than that by WENO3-Z and especially WENO3-JS; however, WENO3-ZM evidently outperforms aforementioned WENO3-Z improvements on resolving the instability along the slip line.

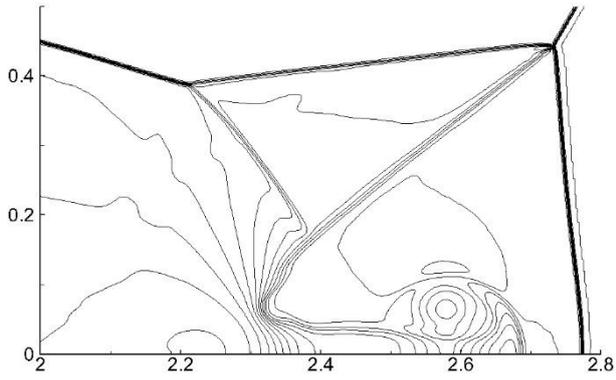
(a) WENO3-JS

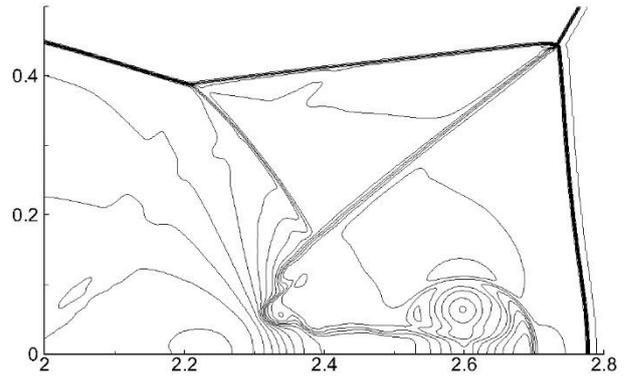
(b) WENO3-Z

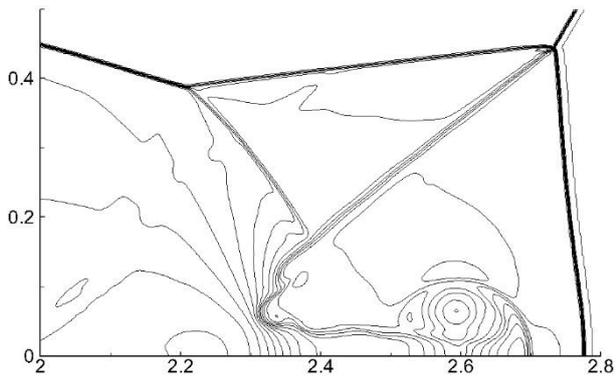
(c) WENO-NP3

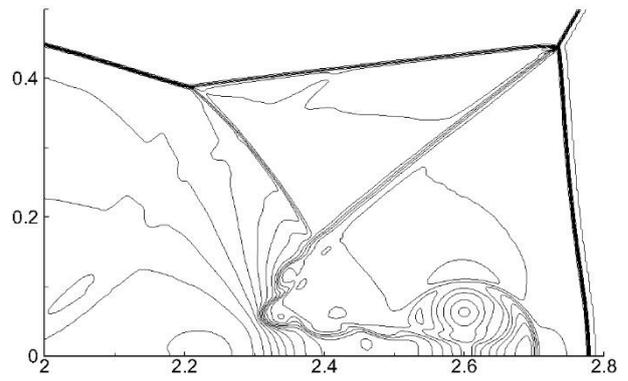
(d) WENO-F3

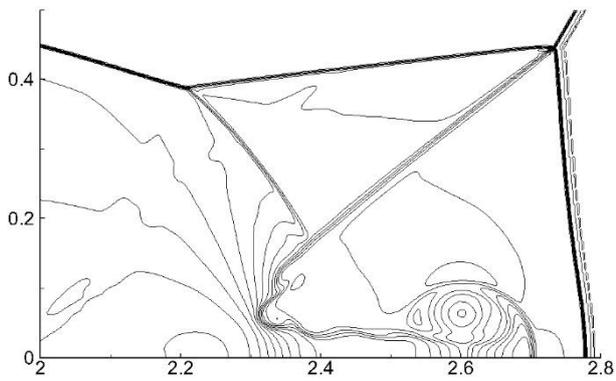
(e) WENO-NN3

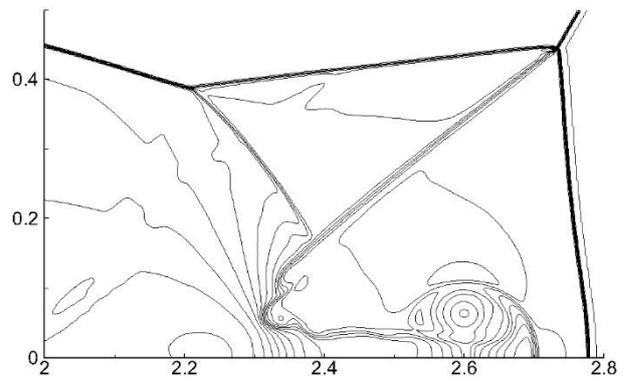
(f) WENO-PZ3

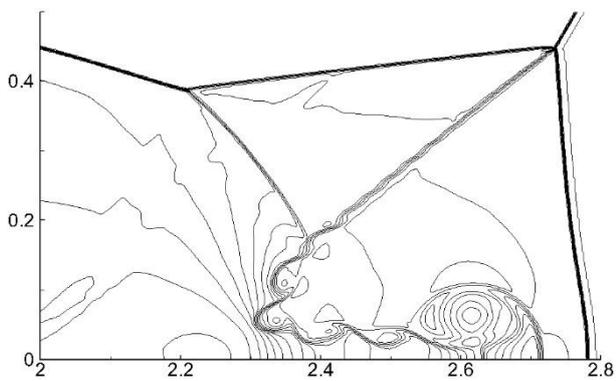
(g) WENO-P+3

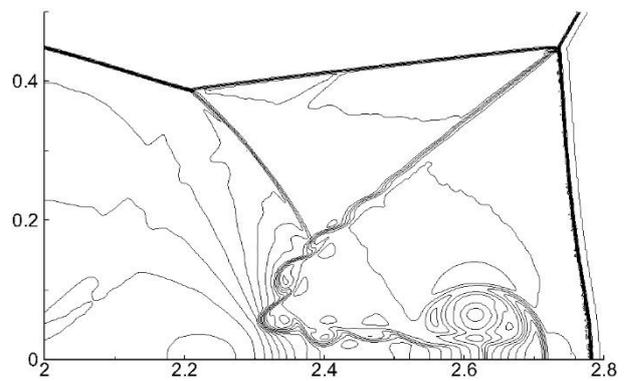
(h) WENO3-ZM

Fig. 10. Density contours of double Mach reflection at on grids $1920 \times 480$ at $t = 2$ by WENO3-ZM and –$Z_{ES}$ with the comparisons of WENO3-JS, WENO3-Z and its other improvements (33 contours from the range of 1.4 to 24)

(3) Time-cost comparison

One may concern the computation cost of proposed schemes with respect to typical third-order WENO-Z schemes. We take use of 2-D Riemann problem and run the computation for 100 steps. The comparative computation costs are: supposing the cost of WENO3-JS is 100, then that of WENO3-Z is 100.05, that of WEN3-ZM is 144.37 and that of WENO3-$Z_{ES}$ is 122.21.

## 6. Conclusions and discussions

In this study, we investigate the accuracy analysis of the third-order WENO-Z scheme in the occurrence of critical points and current WENO3-Z improvements. Based on acquired understandings, two scale-independent, third-order WENO-Z schemes are proposed which could truly preserve the optimal order at critical points. The following conclusions are drawn as:

(1) WENO3-Z improvements referred herein have inherent theoretical deficiency as scale dependency, which might yield undesirable results when extreme variable and/or length-scales are employed. Besides, if a scheme is dependent on physical length-scale, it is difficult to apply it under computational coordinate system.

(2) Although aforementioned WENO3-Z improvements have claimed to recover optimal order at first-order critical points, they actually cannot fulfil the recovery when the points occur at the half nodes.

(3) The prevailing analyses in WENO3-Z improvements have considered the occurrence of critical point and assumed its occurrence at the grid node (such as $x_j$ in Eq. (2)). Such consideration is not sufficient and sometimes would yield incorrect conclusion. Following the former achievements in Ref. [14], a generic analysis is proposed herein which thinks the critical point occurring within grid interval, and the validity of which is testified theoretically and numerically.

(4) Based on achieved theoretical outcomes, two scale-independent, third-order WENO-Z schemes are proposed which can truly recover the optimal order at critical points, namely WENO3-ZM and WENO3-$Z_{ES}$. WENO3-ZM is derived by employing the only solution of $\tau_{CP_1}$ on the stencil $\{x_{j-1}, x_j, x_{j+1}, x_{j+2}\}$ and incorporating with the mapping function. Numerical tests have validated its optimal order recovery at the first-order critical point ($f' = 0, f'' \& f''' \neq 0$); moreover, WENO3-ZM outperforms aforementioned third-order WENO-Z improvements in terms of numerical resolution. WENO3-$Z_{ES}$ is derived by expanding the stencil to that of WENO5-JS, upgrading its smoothness indicators to the two of WENO5-JS and employing a newly derived $\tau_{CP_2}$. Although WENO3-$Z_{ES}$ can achieve the optimal order at even the second-order critical point ($f' = f'' = 0, f''' \neq 0$), it indicates weak robustness in 1-D and 2-D problems by Euler equations and therefore is mainly of theoretical significance.

It is conceivable that the similar problem as above "(3)" exists for other WENO-Z schemes such as WENO5-Z. In this regard we have made studies as well and solutions are acquired. However, due to limited space of the paper, such practices are planned to discuss in other places.

**Acknowledgements**

This study is sponsored by the project of National Numerical Wind-tunnel of China under the grant number NNW2019ZT4-B12.

## Appendix I. Coefficients of WENO3, 5-JS

For reference, the coefficients of candidate schemes, linear weights, and coefficients of smoothness indicators of WENO3,5-JS which correspond to Eqns. (3) -(5) are tabulated in Tables 11-12.

Table 12 Coefficients $a_{kl}^r$ and $d_k$ of candidate schemes $q_k^r$ of WENO-JS with $r=2, 3$

| $r$ | $k$ | $a_{k0}^r$ | $a_{k1}^r$ | $a_{k2}^r$ | $d_k$ |
|---|---|---|---|---|---|
| 2 | 0 | -1/2 | 3/2 | - | 1/3 |
|   | 1 | 1/2 | 1/2 | - | 2/3 |
| 3 | 0 | 2/6 | -7/6 | 11/6 | 1/10 |
|   | 1 | -1/6 | 5/6 | 2/6 | 6/10 |
|   | 2 | 2/6 | 5/6 | -1/6 | 3/10 |

Table 13 Coefficients $b_{kml}^r$ of smoothness indicators of WENO-JS in Eq. (5) with $r=2, 3$

| $r$ | $k$ | $m$ | $b_{km0}^r$ | $b_{km1}^r$ | $b_{km2}^r$ |
|---|---|---|---|---|---|
| 2 | 0 | 0 | -1 | 1 | - |
|   | 1 | 0 | -1 | 1 | - |
| 3 | 0 | 0 | 1 | -4 | 3 |
|   |   | 1 | 1 | -2 | 1 |
|   | 1 | 0 | -1 | 0 | 1 |
|   |   | 1 | 1 | -2 | 1 |

## Appendix II. Proofs of propositions

Suppose $S_C$ and $S_D$ are two substencils of the same pattern, where the subscripts "C, D" indicates the variable distribution at $S_C$ is smoother than at $S_D$, or $\beta_{k,C}^{(r)} < \beta_{k,D}^{(r)}$, and suppose $\alpha_k$ is the non-normalized weight while $\omega_k$ is the normalized weight.

**Proposition 1** Consider $\alpha_k$ with the form $\alpha_k = d_k \left(1 + \tau/\beta_k^{(r)p_i}\right)$, then $\left[\frac{\omega_{k,D}}{\omega_{k,C}}\right]_{p_1} > \left[\frac{\omega_{k,D}}{\omega_{k,C}}\right]_{p_2}$ where $\tau > 0$, $0 < p_1 < p_2 \leq 1$ and $0 < \beta_{k,C}^{(r)} < \beta_{k,D}^{(r)} \leq e$.

**Proof** Observing $\left[\frac{\omega_{k,D}}{\omega_{k,C}}\right]_{p_i} = \left[\frac{\alpha_{k,D}}{\alpha_{k,C}}\right]_{p_i} = \left[1 + \tau\left(\beta_{k,D}^{(r)}\right)^{-p_i}\right] / \left[1 + \tau\left(\beta_{k,C}^{(r)}\right)^{-p_i}\right]$ and $\beta_{k,C}^{(r)} < \beta_{k,D}^{(r)}$, we consider the following function: $y(x) = \frac{1+\tau b^x}{1+\tau a^x}$, where $-1 \leq x < 0$ and $b > a > 0$. It can be seen that $y' = \frac{\tau}{(1+\tau a^x)^2}[(1 + \tau a^x)b^x \ln(b) - (1 + \tau)a^x \ln(a)]$, where $\tau/(1 + \tau a^x)^2 > 0$. Let $g(m) = \frac{m^x \ln(m)}{1+\tau m^x}$, then one can see that $g(b) > g(a)$ indicates $[(1 + \tau a^x)b^x \ln(b) - (1 + \tau b^x)a^x \ln(a)] > 0$. Therefore, the derivative of $g(m)$ should be checked, which can be derived as $g' = \frac{x \ln(m)+1+\tau m^x}{m(1+\tau m^x)^2 m^{-x}}$. one can see that the numerator of $g'$, $x \ln(m) + 1 + \tau m^x$, is a monotonously decreasing function of $m$ for $x \in [-1,0)$. Because $m^x > 0$, $m_0$ which is the root of the numerator being zero will increase with the increase of $\tau$. The minimum of $m_0$ will be attained when $\tau = 0$ as $m_0 = e^{-\frac{1}{x}}$. Considering $-1 \leq x < 0$, then $\min(m_0) = e^{-\frac{1}{x}}|_{-1} = e$. Hence, $g' > 0$ holds when $0 < m < e$ and $-1 \leq x < 0$. Moreover, if $0 < x = a < x = b \leq e$, $g' > 0$ indicates $g(b) > g(a)$ and thereby $y' > 0$. Consequently, $\left[\frac{\omega_{k,D}}{\omega_{k,C}}\right]_{p_1} - \left[\frac{\omega_{k,D}}{\omega_{k,C}}\right]_{p_2} = y(-p_1) - y(-p_2)|_{-p_1 > -p_2} > 0$, or $\left[\frac{\omega_{k,D}}{\omega_{k,C}}\right]_{p_1} > \left[\frac{\omega_{k,D}}{\omega_{k,C}}\right]_{p_2}$.

**Proposition 2.** Consider $\alpha_k$ with the form $\alpha_k = d_k\left(1 + \tau^{p_i}/\beta_k^{(r)}\right)$, then $\left[\frac{\omega_{k,D}}{\omega_{k,C}}\right]_{p_1} < \left[\frac{\omega_{k,D}}{\omega_{k,C}}\right]_{p_2}$ for $1 \leq p_1 < p_2$ and $0 < \tau < 1$; $\left[\frac{\omega_{k,D}}{\omega_{k,C}}\right]_{p_1} > \left[\frac{\omega_{k,D}}{\omega_{k,C}}\right]_{p_2}$ for $1 \leq p_1 < p_2$ and $1 < \tau$; $\left[\frac{\omega_{k,D}}{\omega_{k,C}}\right]_{p_1} = \left[\frac{\omega_{k,D}}{\omega_{k,C}}\right]_{p_2}$ for $1 \leq p_1 < p_2$ and $\tau = 1$.

**Proof** Observing $\left[\frac{\omega_{k,D}}{\omega_{k,C}}\right]_{p_i} = \left[\frac{\alpha_{k,D}}{\alpha_{k,C}}\right]_{p_i} = \left[1 + \tau^{p_i}/\beta_{k,D}^{(r)}\right]/\left[1 + \tau^{p_i}/\beta_{k,C}^{(r)}\right]$ and $\beta_{k,C}^{(r)} < \beta_{k,D}^{(r)}$, we consider the following function: $y(x) = \frac{1+\tau^x/b}{1+\tau^x/a}$ where $x \geq 1$ and $b > a > 0$. It can be seen that $y' = \frac{a\left(1+\frac{\tau^x}{b}\right)t^x}{(b+\tau^x)(a+\tau^x)^2}\ln(\tau)(a-b)$, where $\frac{a\left(1+\frac{\tau^x}{b}\right)t^x}{(b+\tau^x)(a+\tau^x)^2} > 0$. Then

(1) If $0 < \tau < 1$, then $y' > 0$, $\left[\frac{\omega_{k,D}}{\omega_{k,C}}\right]_{p_1} - \left[\frac{\omega_{k,D}}{\omega_{k,C}}\right]_{p_2} = \frac{1+\frac{\tau^{p_1}}{\beta_{k,D}^{(r)}}}{1+\frac{\tau^{p_1}}{\beta_{k,C}^{(r)}}} - \frac{1+\frac{\tau^{p_2}}{\beta_{k,D}^{(r)}}}{1+\frac{\tau^{p_2}}{\beta_{k,C}^{(r)}}} = y(p_1) - y(p_2) < 0$.

(2) If $1 < \tau$, then $y' < 0$, $\left[\frac{\omega_{k,D}}{\omega_{k,C}}\right]_{p_1} - \left[\frac{\omega_{k,D}}{\omega_{k,C}}\right]_{p_2} = \frac{1+\frac{\tau^{p_1}}{\beta_{k,D}^{(r)}}}{1+\frac{\tau^{p_1}}{\beta_{k,C}^{(r)}}} - \frac{1+\frac{\tau^{p_2}}{\beta_{k,D}^{(r)}}}{1+\frac{\tau^{p_2}}{\beta_{k,C}^{(r)}}} = y(p_1) - y(p_2) > 0$.

(3) If $\tau = 1$, then $y' = 0$, $\left[\frac{\omega_{k,D}}{\omega_{k,C}}\right]_{p_1} - \left[\frac{\omega_{k,D}}{\omega_{k,C}}\right]_{p_2} = \frac{1+\frac{\tau^{p_1}}{\beta_{k,D}^{(r)}}}{1+\frac{\tau^{p_1}}{\beta_{k,C}^{(r)}}} - \frac{1+\frac{\tau^{p_2}}{\beta_{k,D}^{(r)}}}{1+\frac{\tau^{p_2}}{\beta_{k,C}^{(r)}}} = y(p_1) - y(p_2) = 0$.

**Proposition 3.** Consider $\alpha_k$ with the form $\alpha_k = d_k\left(1 + c_i \times \tau/\beta_k^{(r)}\right)$, then $\left[\frac{\omega_{k,D}}{\omega_{k,C}}\right]_{c_1} < \left[\frac{\omega_{k,D}}{\omega_{k,C}}\right]_{c_2}$ for $c_1 > c_2$.

**Proof** Observing $\left[\frac{\omega_{k,D}}{\omega_{k,C}}\right]_{ci} = \left[\frac{\alpha_{k,D}}{\alpha_{k,C}}\right]_{ci} = \frac{1+c_i\tau/\beta_{k,D}^{(r)}}{1+c_i\tau/\beta_{k,C}^{(r)}}$, we consider the function: $y(x) = \frac{1+x\cdot\tau/b}{1+x\cdot\tau/a}$ where $b > a > 0$. It can be seen that $y' = \frac{\tau}{ab(1+\tau x/a)^2}(a-b) < 0$. Hence if $x = c_1 > x = c_2$, it holds that $\left[\frac{\omega_{k,D}}{\omega_{k,C}}\right]_{c_1} - \left[\frac{\omega_{k,D}}{\omega_{k,C}}\right]_{c_2} = \frac{1+c_1\frac{\tau}{b}}{1+c_1\frac{\tau}{a}} - \frac{1+c_2\frac{\tau}{b}}{1+c_2\frac{\tau}{a}} = y(c_1) - y(c_2) < 0$ considering $y' < 0$.

**Proposition 4.** Consider $\alpha_k$ with the form $\alpha_k = d_k\left(1 + \left(\tau/\beta_k^{(r)}\right)^{p_i}\right)$, then $\left[\frac{\omega_{k,D}}{\omega_{k,C}}\right]_{p_1} < \left[\frac{\omega_{k,D}}{\omega_{k,C}}\right]_{p_2}$ for $p_1 > p_2$ providing $0.278 < \frac{\tau}{\beta_{k,D}^{(r)}} < \frac{\tau}{\beta_{k,C}^{(r)}}$.

**Proof** Observing $\left[\frac{\omega_{k,D}}{\omega_{k,C}}\right]_{p_i} = \left[\frac{\alpha_{k,D}}{\alpha_{k,C}}\right]_{p_i} = \frac{1+\left(\tau/\beta_{k,D}^{(r)}\right)^{p_i}}{1+\left(\tau/\beta_{k,C}^{(r)}\right)^{p_i}}$ where $\tau/\beta_{k,D}^{(r)} < \tau/\beta_{k,C}^{(r)}$, we consider the following function: $y(x) = \frac{1+a^x}{1+b^x}$ where $b > a$. Accordingly, $y' = a^x \ln(a)(1+b^x) - (1+a^x)b^x \ln(b)$. Supposing $g(m) = \frac{m^x \ln(m)}{1+m^x}$, then $g(a) > g(b)$ indicates $y' = [a^x \ln(a)(1+b^x) - (1+a^x)b^x \ln(b)] > 0$. As derived, $g' = \frac{x\ln(m)+1+m^x}{m(1+m^x)^2 m^{-x}}$, therefore if $g' > 0$, $g(b) > g(a)$ providing $b > a$; and vice vesa. Moreover, the numerator of $g'$, namely $x \ln(m) + 1 + m^x$, can be found to monotonously increase with $m$, and the root of the zero point is $m_0 = e^{-\frac{1.2784645428}{x}}$.

Hence $m_0$ increase with that of $x$, and $\max(m_0) = e^{-\frac{1.2784645428}{x}}|_{x=1} = 0.278$ thereby. So when $m \geq 0.278$, $g' \geq 0$ and $y' \leq 0$. In short, when $0.278 < \frac{\tau}{\beta_{k,D}^{(r)}} < \frac{\tau}{\beta_{k,C}^{(r)}}$, $y(p_1) < y(p_2)$, or $\left[\frac{\omega_{k,D}}{\omega_{k,C}}\right]_{p_1} < \left[\frac{\omega_{k,D}}{\omega_{k,C}}\right]_{p_2}$.

**Proposition 5.** Consider the generic quadratic form of $f$ as $\tau(f) = (f_{j-1}, f_j, f_{j+1})[a_{i_1 i_2}](f_{j-1}, f_j, f_{j+1})^T$ where $[a_{i_1 i_2}]$ is a $3 \times 3$ matrix and $i_1, i_2 = 1...3$ Supposing the first-order critical point would occur at $x_c = x_j + \lambda \cdot \Delta x$ where $-1 < \lambda < 1$ and $\{f'_{x_c} = 0, f''_{x_c} \& f'''_{x_c} \neq 0\}$, then no solution of $a_{i_1 i_2}$ exists such that Taylor expansion of $\tau(f)$ toward $x_c$ has the leading error as $O(\Delta x^5)$.

**Proof** Expand $\tau(f)$ toward $x_c$ by Taylor expansion and obtain: $\tau(f) \approx \sum_{n=0}^{\infty} \Delta x^n g(n, f, \lambda, a_{i_1 i_2})$ where $g(n, f, \lambda, a_{i_1 i_2})$ is the function of $n, f, \lambda$ and $a_{i_1 i_2}$. Given the leading error order of $\tau(f)$ as $k$ and the condition $\{f'_{x_c} = 0, f''_{x_c} \& f'''_{x_c} \neq 0\}$, if the following equations have the solutions about $a_{i_1 i_2}$

$$\begin{cases} g(0, f, \lambda, a_{i_1 i_2}) = 0 \\ g(1, f, \lambda, a_{i_1 i_2}) = 0 \\ g(2, f, \lambda, a_{i_1 i_2}) = 0 \\ \cdots \\ g(k-1, f, \lambda, a_{i_1 i_2}) = 0 \end{cases},$$

Then $\tau(f) = O(\Delta x^k)$ will hold. Take $g(2, f, \lambda, a_{i_1 i_2})$ as an example which can be derived as

$g(2, f, \lambda, a_{i_1 i_2}) = \left[(\sum_{i_1=1}^{3} \sum_{i_2=1}^{3} a_{i_1 i_2}) \lambda^2 + (-2a_{33} + a_{21} - a_{23} + 2a_{11} + a_{12} - a_{32})\lambda + \frac{1}{2}a_{23} + a_{11} + a_{31} + \frac{1}{2}a_{12} + +\frac{1}{2}a_{32} + a_{13} + a_{33} + \frac{1}{2}a_{21}\right] f_{j-\lambda} f''_{j-\lambda}$

where $\{f'_{x_c} = 0, f''_{x_c} \& f'''_{x_c} \neq 0\}$ has been casted. In order to get the solution of $a_{i_1 i_2}$ such that $g(2, f, \lambda, a_{i_1 i_2}) = 0$, the coefficients regarding $\lambda$ should be zero, or the following equations should be satisfied:

$$\begin{cases} \sum_{i_1=1}^{3} \sum_{i_2=1}^{3} a_{i_1 i_2} = 0 \\ -2a_{33} + a_{21} - a_{23} + 2a_{11} + a_{12} - a_{32} = 0 \\ \frac{1}{2}a_{23} + a_{11} + a_{31} + \frac{1}{2}a_{12} + +\frac{1}{2}a_{32} + a_{13} + a_{33} + \frac{1}{2}a_{21} = 0 \end{cases}$$

Taking $k$ as 5, the whole derived equations can be defined and solved, and $\{a_{i_1 i_2}\}$ can be acquired thereby. Substituting the obtained $\{a_{i_1 i_2}\}$ into $\tau(f)$, it is found $\tau(f) = 0$ or no meaningful solution of $\{a_{i_1 i_2}\}$ exists.

**Proposition 6.** Consider the generic quadratic form of $f$ as $\tau(f) = (f_{j-1}, f_j, f_{j+1}, f_{j+2})[a_{i_1 i_2}](f_{j-1}, f_j, f_{j+1}, f_{j+2})^T$ where $[a_{i_1 i_2}]$ is a $4 \times 4$ matrix and $i_1, i_2 = 1\ldots 4$. Supposing the first-order critical point would occur at $x_c = x_j + \lambda \cdot \Delta x$ where $-1 < \lambda < 2$ and $\{f'_{x_c} = 0, f''_{x_c} \& f'''_{x_c} \neq 0\}$, the solution of $\tau(f)$ is wanted such that its Taylor expansion toward $x_c$ having the leading error as: (1) $O(\Delta x^5)$ when $\lambda \neq 0 \& -\frac{1}{2}$; (2) $O(\Delta x^7)$ when $\lambda = -\frac{1}{2}$, then the only solution is found as: $\tau(f) = c \times (-f_{j+2} + 3f_{j+1} + 21f_j - 23f_{j-1}) \times (f_{j+2} - 3f_{j+1} + 3f_j - f_{j-1})$ where $c$ is a free parameter.

**Proof** Expand $\tau(f)$ toward $x_c$ by Taylor expansion and obtain: $\tau(f) \approx \sum_{n=0}^{\infty} \Delta x^n g(n, f, \lambda, a_{i_1 i_2})$ where $g(n, f, \lambda, a_{i_1 i_2})$ is the function of $n, f, \lambda$ and $a_{i_1 i_2}$. The following equations are solved first to derive $\{a_{i_1 i_2}\}$ such that the leading error of $\tau(f)$ would be $O(\Delta x^5)$ under $\{f'_{x_c} = 0, f''_{x_c} \& f'''_{x_c} \neq 0\}$

$$\begin{cases} g(0, f, \lambda, a_{i_1 i_2}) = 0 \\ g(1, f, \lambda, a_{i_1 i_2}) = 0 \\ g(2, f, \lambda, a_{i_1 i_2}) = 0 \\ g(3, f, \lambda, a_{i_1 i_2}) = 0 \\ g(4, f, \lambda, a_{i_1 i_2}) = 0 \end{cases}$$

where $\lambda$ is regarded as arbitrary parameter. The solving processes are similar and referred to that in Proposition 5. Next, by means of substituting $\{a_{i_1 i_2}\}$ into $\tau(f)$ and expanding it toward a specific $x_c$ where $\lambda = -\frac{1}{2}$: $x_c = x_j - \frac{1}{2} \cdot \Delta x$, $g\left(5, f, \lambda = -\frac{1}{2}, a_{i_1 i_2}\right)$ and $g\left(6, f, \lambda = -\frac{1}{2}, a_{i_1 i_2}\right)$ can be derived. Then, the following equations are solved under $f'_{j-1/2} = 0$

$$\begin{cases} g\left(5, f, \lambda = -\frac{1}{2}, a_{i_1 i_2}\right) = 0 \\ g\left(6, f, \lambda = -\frac{1}{2}, a_{i_1 i_2}\right) = 0 \end{cases},$$

The previous solution $\{a_{i_1 i_2}\}$ is further specified as

$$[a_{i_1 i_2}] = \begin{bmatrix} -23a_{44} & 90a_{44} - a_{21} & -66a_{44} - a_{31} & 22a_{44} - a_{41} \\ a_{21} & -63a_{44} & 54a_{44} - a_{32} & -a_{24} - 18a_{44} \\ a_{31} & a_{32} & 9a_{44} & -a_{43} - 6a_{44} \\ a_{41} & a_{42} & a_{43} & a_{44} \end{bmatrix}$$

Considering the equivalence property of quadratic form, the above $[a_{i_1 i_2}]$ is equivalent to:

$$[a_{i_1 i_2}] = \begin{bmatrix} -23a_{44} & 90a_{44} & -66a_{44} & 22a_{44} \\ 0 & -63a_{44} & 54a_{44} & -18a_{44} \\ 0 & 0 & 9a_{44} & -6a_{44} \\ 0 & 0 & 0 & a_{44} \end{bmatrix} = a_{44} \begin{bmatrix} -23 & 90 & -66 & 22 \\ 0 & -63 & 54 & -18 \\ 0 & 0 & 9 & -6 \\ 0 & 0 & 0 & 1 \end{bmatrix}$$

By the substitution of $[a_{i_1 i_2}]$ and after factorization, $\tau(f)$ can be finally derived as:
$\tau(f) = -a_{44} \times (-f_{j+2} + 3f_{j+1} + 21f_j - 23f_{j-1}) \times (f_{j+2} - 3f_{j+1} + 3f_j - f_{j-1})$